\documentclass[a4paper,12pt]{amsart}

\usepackage{amsmath,amsfonts,amssymb,amsthm}
\usepackage{graphics}
\usepackage{epsfig}
\usepackage{esint}
\usepackage{hyperref}

\newcommand{\ignore}[1]{}

\newtheorem{proposition}{Proposition}
\newtheorem{theorem}{Theorem}

\newtheorem{lemma}{Lemma}

\usepackage{color}
\definecolor{darkred}{rgb}{0.9,0.1,0.1}
\definecolor{darkblue}{rgb}{0,0,0.7}
\definecolor{darkgreen}{rgb}{0,0.5,0}

\begin{document}

\title[Thresholding and minimizing movements]{The thresholding scheme for mean curvature flow and
de Giorgi's ideas for minimizing movements}
\author{Tim Laux}
\author{Felix Otto}

\begin{abstract}
	We consider the thresholding scheme and explore its connection to De Giorgi's ideas on gradient flows in metric spaces; here applied to mean curvature flow as the steepest descent of the interfacial area.
 The basis of our analysis is the observation by Esedo\u{g}lu and the second author that thresholding can be interpreted as a minimizing movements scheme for an energy that approximates the interfacial area. 
 De Giorgi's framework provides an optimal energy dissipation relation for the scheme in which we pass to the limit to derive a dissipation-based weak formulation of mean curvature flow.
 Although applicable in the general setting of arbitrary networks, here we restrict ourselves to the case of  a single interface, which allows for a compact, self-contained presentation.  
\end{abstract}
\maketitle
\section{Introduction and context}

The purpose of these notes is to draw a connection between De Giorgi's tools for
minimizing movements, that is, gradient flows in metric spaces on the one
hand, and the very popular thresholding scheme for flow of a hyper-surface by its mean curvature
on the other hand.
While we have developed this connection in the case of multiple phases
with surface energies and mobilities depending on the pair of phases,
as is relevant for grain growth in polycrystals, and when the notion
of viscosity solution is not available, we present our results
here in the simplest setting of two phases. Our
presentation is essentially self-contained.

\smallskip

What makes the evolution of the boundary $\partial\Omega$ of a set $\Omega$
by its mean curvature $H$ valuable for modeling in materials science
is that it is driven by the reduction of the (total) interfacial area of $\partial\Omega$,
which relies on the mean curvature $H$, the sum of the principal
curvatures, being the first variation of the interfacial area.
There is a more intimate connection between mean curvature flow (MCF) and the functional $E$ of interfacial area
of a configuration: MCF can formally be understood as a gradient flow of $E$. 
We stress that a dynamical system that can be written as a gradient flow, that is,
a steepest descent in an energy landscape, does not just rely on the height function $E$,
but also on a notion of distance on configuration space, which is typically described by
a metric tensor $g$ in the sense of Riemannian geometry. In case of MCF, the tangent
space in some configuration $\Omega$ should be thought of as consisting of all normal
velocities $V$, i.e.,  functions on $\partial \Omega$, while the configuration-dependent metric tensor $g_\Omega$ is
given by the $L^2$-inner product on $\partial\Omega$. 

\smallskip

Still formally, any gradient flow allows for a natural discretization in time.
Every step of the discretization comes in form of a variational problem, just involving
the functional $E$ and the induced distance $d$, cf.~(\ref{wg66}), but not the metric tensor $g$ and
the differential of $E$ -- it thus relies rather on the ``metric'', but not the
differential structure. Following De Giorgi, we call such a scheme a minimizing
movements scheme. We recall that, as in elementary differential geometry,
the induced distance $d$ on a Riemannian manifold $({\mathcal M},g)$ is defined
via $d^2(\chi_0,\chi_1)$ $:=\inf\{\int_0^1 g_{\chi_s}(\frac{d\chi}{ds},\frac{d\chi}{ds})ds\}$,
where the infimum is taken over all curves $[0,1]\ni s\mapsto\chi_s$ connecting
$\chi_0$ to $\chi_1$ (we use the letter $\chi$ because we think of a characteristic
function describing the configuration). 
We note that in the Euclidean case, the Euler-Lagrange equation
of (\ref{wg66}) turns into the implicit Euler scheme for $\frac{d\chi}{dt}=\mathop{\rm grad}E|_{\chi}$.

\smallskip

However, this infinite-dimensional Riemannian structure making MCF 
a gradient flow leads to a degenerate induced metric (i.e.~$d\equiv0$): It can be seen that the 
infimum of $\int_0^1\int_{\partial\Omega_s}V_s^2 ds$
over all curves of configurations $[0,1]\ni s\mapsto\Omega_s$ with normal velocity 
$[0,1]\ni s\mapsto V_s$,
connecting some given configurations $\Omega_0$ and $\Omega_1$, vanishes \cite{MichorMumford}.

\smallskip

Nonetheless, a minimizing movements scheme (before the latter) for MCF has
been formulated by Almgren~et.~al.~\cite{AlmgrenTaylorWang}, with $E(\Omega)$ being
the surface area of $\partial\Omega$ and with $d^2(\Omega_1,\Omega_0)$ 
$=4\int_{\Omega_1\triangle\Omega_0}{\rm dist}(\cdot,\partial\Omega_0)$. 
Luckhaus~et.~al.~\cite{LuckhausSturzenhecker} have established a (long-time) convergence
result for this scheme. This convergence result is conditional
in the sense that a condition like in (\ref{wg04}) has to be imposed.

\smallskip

Thresholding, cf.~(\ref{wg59}), is a very well performing and widely used numerical scheme for MCF,
introduced by Osher~et.~al.~\cite{BenceMerrimanOsher}. Also the convolution step,
which after spatial discretization can be carried out by the Fast Fourier Transform,
is of low complexity. Right from the beginning, thresholding
has attracted the attention of analysts; since it obviously conserves the comparison principle
for MCF, it has been shown to converge to MCF in the sense of viscosity solution 
in the two-phase case \cite{Evans}.

\smallskip

Esedo\u{g}lu and the second author \cite{EsedogluOtto} realized that thresholding also respects the 
gradient-flow structure of MCF, in the sense that it can be interpreted as a minimizing
movements scheme, cf.~Lemma \ref{Le2}. This was used in the multi-phase case
to extend thresholding to surface tensions and mobilities \cite{EsedogluSalvador}
that depend on the pair of grains, while keeping its low complexity. 
It was also used by the present authors to 
provide several types of convergence results; presently, all of them are
conditional in the sense of assumption (\ref{wg04}), in the tradition of 
\cite{LuckhausSturzenhecker}. 

\smallskip

The first result \cite{LauxOtto1} provided the same
limiting notion of solution for MCF as in \cite{LuckhausSturzenhecker}. However, this
weak notion of solution does not imply the dissipation inequality natural to a gradient flow.
It is Brakke's weak notion of solution for MCF that is based on a localization of
the dissipation inequality; in \cite{LauxOtto2}, we establish a (still conditional)
convergence result towards this inequality-based notion of solution.

\smallskip

For any gradient flow in a Riemannian context $({\mathcal M},g,E)$, 
there is yet another notion of weak solution based on a single inequality,
namely $E(\chi(T))$ $+\int_0^T\frac{1}{2}g_{\chi}(\frac{d\chi}{dt},\frac{d\chi}{dt})
+\frac{1}{2}|{\rm grad} E_{|\chi}|^2dt$ $\le E(\chi(0))$. This elementary
observation is credited to De Giorgi; its appeal lies in the fact that
it is potentially more stable in limiting procedures because only lower semi-continuity
is needed (as provided by Propositions \ref{Pr2} and \ref{Pr3}). 
The main result of this paper, Theorem \ref{Th}, precisely establishes
this inequality in the case of MCF, cf.~(\ref{wg62}).

\smallskip

One advantage of a minimizing movements scheme, cf.~(\ref{wg66}), lies in the fact that it
automatically comes with the a priori estimate 
$E(\chi^N)+\sum_{n=1}^N\frac{1}{2h}d^2(\chi^n,\chi^{n-1})$ $\le E(\chi^0)$, which
is obtained by using $\chi^{n-1}$ as a competitor in (\ref{wg66}). In the limit $h\downarrow0$, 
this inequality formally turns into 
$E(\chi(T))$ $+\int_0^T\frac{1}{2}g_{\chi}(\frac{d\chi}{dt},\frac{d\chi}{dt})dt$ $\le E(\chi(0))$,
which misses the formally correct identity by a factor of $2$. On the level of the metric
structure, De Giorgi provides tools to capture the missing term 
$\int_0^T\frac{1}{2}|{\rm grad} E_{|\chi}|^2dt$, see Lemma \ref{Le1}. We take
the proof from the monograph \cite{AmbrosioGigliSavare}. 

\smallskip

As a consequence of these notions and tools of De Giorgi, our (conditional) convergence proof
for the thresholding scheme in fact is rather ``soft'', softer than \cite{LauxOtto1} which
relied on the notion of tilt excess and the fine structure of Caccioppoli sets, 
and certainly softer than \cite{LuckhausSturzenhecker}
which relied on regularity theory for minimal surfaces. We believe that these tools
have a wider potential for geometric evolutions or non-linear PDE of gradient-flow type.
For the broader context and more references, we refer to \cite{LauxOtto1}.

\section{Main result and structure of proof}

Given an initial configuration, as described by its (Lebesgue-\-mea\-sur\-able)
characteristic function $\chi^0\colon\mathbb{R}^d\rightarrow\{0,1\}$,
and a time step size $h>0$,
the thresholding scheme iteratively produces configurations at time steps $n=1,2,\cdots$,
encoded by their characteristic functions $\chi^n$, via convolution and ``thresholding'':
\begin{align}\label{wg59}
\chi^n:=\left\{\begin{array}{cl}
1&\mbox{where}\;G_h*\chi^{n-1}>\frac{1}{2}\\
0&\mbox{else}\end{array}\right\},
\end{align}
where $G_h$ denotes the heat kernel at time $\frac{h}{2}$, that is,
\begin{align}\label{wg76}
G_h(z):=\frac{1}{\sqrt{2\pi h}^d}\exp(-\frac{|z|^2}{2h})
\end{align}
(like in stochastic analysis, we take $\frac{h}{2}$ so that $G_1$ is the standard Gaussian).
We interpolate piecewise constant in time:
\begin{align}\label{wg61}
\chi_h(t)=\chi^n\quad\mbox{for}\;t\in[nh,(n+1)h),
\quad\chi_h(t)=\chi^0\quad\mbox{for}\;t\le 0.
\end{align}
For simplicity, we pass from the whole space $\mathbb{R}^d$ to a torus as the spatial domain; 
by rescaling, we may w.~l.~o.~g.~take the unit torus $[0,1)^d$.
We also restrict to the finite time horizon $T<\infty$ and $h\le 1$.

\medskip

Our main result is the following convergence result, which is only a conditional one
since assumption (\ref{wg04}) on the energies $E_h$ defined below in (\ref{wg06}) presumably cannot be verified. It is the opposite
direction, 
$c_0\int_{(0,T)\times[0,1)^d}|\nabla\chi|dt$ 
$\le\liminf_{h\downarrow0}$ $\int_0^T E_h(\chi_h(t))dt$,
that follows from (\ref{wg07}), see (\ref{wg84}) in Lemma \ref{Le3}. 
Here and in the sequel, $|\nabla\chi|dt$ denotes the total variation of 
the distribution $\nabla\chi$ in $(0,T)\times[0,1)^d$, provided the latter is a bounded measure.
This notation is justified since in the present case, $|\nabla\chi|dt$ is
(Lebesgue) equi-integrable in $t$ and thus admits a density $|\nabla\chi|$.
In the sequel, $\nu\in L^1(|\nabla\chi|dt)$ denotes the measure-theoretic normal
(characterized through the polar factorization $\nabla\chi=\nu|\nabla\chi|dt$
and $|\nu|=1$ $|\nabla\chi|dt$-almost everywhere).

\begin{theorem}\label{Th}
Given $\chi^0$ as above and such that $\nabla\chi^0$ is a bounded measure, 
and a sequence $h\downarrow 0$; let $\chi_h$ be defined by (\ref{wg59}) and (\ref{wg61}).
Suppose that there exists a 
$\chi\colon(0,T)\times[0,1)^d\rightarrow[0,1]$ such that
%
\begin{align}
\chi_h\rightharpoonup\chi\quad\mbox{in}\; L^1((0,T)\times[0,1)^d).\label{wg07}
\end{align}
Then we have $\chi\in\{0,1\}$ (Lebesgue)-a.~e. and $\nabla\chi$ is a bounded measure
which is equi-integrable in $t$.
If we assume in addition
\begin{align}
\limsup_{h\downarrow 0}\int_0^T E_h(\chi_h(t))dt\le 
c_0\int_{(0,T)\times[0,1)^d}|\nabla\chi|dt,\label{wg04}
\end{align}
where $c_0:=\frac{1}{\sqrt{2\pi}}$, 
then there exists $H\in L^2(|\nabla\chi|dt)$ with
%
\begin{align}\label{wg60}
	\int_{(0,T)\times[0,1)^d}(\nabla\cdot\xi-\nu\cdot\nabla\xi\nu)|\nabla\chi|dt
	=-\int_{(0,T)\times[0,1)^d} H\nu\cdot\xi|\nabla\chi|dt
\end{align}
for all $\xi\in C^\infty_0((0,T)\times[0,1)^d)^d$,
and $V$ $\in L^2$ $(|\nabla\chi|dt)$ with
%
\begin{align}\label{wg57bis}
\int_{[0,1)^d}\zeta(t=0)\chi^0dx&+\int_{(0,T)\times[0,1)^d}\partial_t\zeta\chi dxdt\nonumber\\
	&+\int_{(0,T)\times[0,1)^d}\zeta V|\nabla\chi|dt=0
\end{align}
%
for all $\zeta\in C^\infty_0([0,T)\times[0,1)^d)$, such that
\begin{align}\label{wg62}
\limsup_{\tau\downarrow 0}&\frac{1}{\tau}
\int_{(T-\tau,T)\times[0,1)^d}|\nabla\chi|dt\nonumber\\
&+\int_{(0,T)\times[0,1)^d}\big(V^2+(\frac{H}{2})^2\big)|\nabla\chi|dt
\le\int_{[0,1)^d}|\nabla\chi^0|.
\end{align}
\end{theorem}

%
We note that in case of $\{\chi=1\}$ being smooth in time-space $[0,T]\times[0,1)^d$,
$|\nabla\chi|$ coincides with the surface measure, $\nu$ with the (inner) normal, and
$H$ and $V$ coincide with the mean curvature (with the convention
that convex sets have positively curved boundary) and normal velocity
(with the convention that growing sets have positive velocity), respectively. In
addition, (\ref{wg57bis}) yields that $\chi(t=0)=\chi^0$. Moreover,
expanding the square $V^2+(\frac{H}{2})^2$ $=(V+\frac{H}{2})^2-VH$,
and appealing to the classical formula $\frac{d}{dt}\int_{[0,1)^d}|\nabla\chi|$ 
$=\int_{[0,1)^d}VH|\nabla\chi|$ (which relies on the fact that mean curvature 
describes the first variation of the surface area), we see that (\ref{wg62}) turns
into $\int_{(0,T)\times[0,1)^d}(V+\frac{H}{2})^2|\nabla\chi|dt\le 0$, and
thus MCF in form of $V=-\frac{H}{2}$ (the factor $\frac{1}{2}$ stems from
the normalization in (\ref{wg76})). Therefore, the inequality (\ref{wg62}) may
be considered a weak notion of MCF.

\medskip

In the sequel, we omit writing the time-space domain $(0,T)\times[0,1)^d$
when integrating the Lebesgue measure $dxdt$ or the limiting surface measure $|\nabla\chi|dt$.
However, the convolution $*$, for which we reserve the $z$-variable, is always w.~r.~t.~
$\mathbb{R}^d$.

\medskip

The next elementary lemma provides the necessary notions and results
on abstract minimizing movements schemes.

\begin{lemma}\label{Le1}
Let $({\mathcal M},d)$ be a compact metric space and $E\colon{\mathcal M}\rightarrow\mathbb{R}$
be continuous. Given $\chi^0\in{\mathcal M}$ and $h>0$ consider a sequence
$\{\chi^n\}_{n\in\mathbb{N}}$ satisfying
\begin{align}\label{wg66}
\chi^n\quad\mbox{minimizes}\quad\frac{1}{2h}d^2(u,\chi^{n-1})+E(u)\quad
\mbox{among all}\;u\in{\mathcal M}.
\end{align}
Then we have for all $t\in\mathbb{N}h$
\begin{align}\label{wg67}
	\lefteqn{E(\chi(t))}\nonumber\\
	&+\frac{1}{2}\int_0^t\big(\frac{1}{h^2}d^2(\chi(s+h),\chi(s))+|\partial E(u(s))|^2\big)ds
\le E(\chi^0).
\end{align}
Here $\chi(t)$ is the piecewise constant interpolation, cf.~(\ref{wg61}),
$u(t)$ is another interpolation (the ``variational interpolation'') satisfying
\begin{align}
\int_0^\infty\frac{1}{2h^2}d^2(u(t),\chi(t))dt&\le E(\chi^0)\label{wg69},\\
E(u(t))&\le E(\chi(t))\;\;\mbox{for all}\;t\ge 0,\label{wg70}
\end{align}
and $|\partial E(u)|$ is the ``metric slope'' defined through
\begin{align}\label{wg26}
|\partial E(u)|:=\limsup_{v:d(v,u)\rightarrow 0}\frac{(E(u)-E(v))_+}{d(v,u)}\in[0,\infty].
\end{align}
\end{lemma}


The next elementary but crucial lemma establishes that the thresholding scheme
is a minimizing movements scheme.

\begin{lemma}\label{Le2}
Expression (\ref{wg59}) satisfies (\ref{wg66}) provided we define
\begin{align}
	{\mathcal M}&:=\{u\colon[0,1)^d\rightarrow[0,1]\;\mbox{measurable}\},\label{wg92}\\
E_h(u)&:=\frac{1}{\sqrt{h}}\int(1-u)\,G_h*u dx,\label{wg06}\\
d_h(u,u')&:=\big(2\sqrt{h}\int|G_\frac{h}{2}*(u-u')|^2dx\big)^\frac{1}{2}.\label{wg68}
\end{align}
Furthermore, $({\mathcal M},d_h)$ is a compact metric space and $E_h$ continuous.
\end{lemma}

We will mostly use (\ref{wg68}) in form of
\begin{align}\label{wg18}
\frac{1}{2h}d_h^2(u,u')=\frac{1}{\sqrt{h}}\int|G_\frac{h}{2}*(u-u')|^2dx.
\end{align}

\medskip


The first part of the next lemma provides compactness.
The second part contains the (only) way we use the convergence assumption (\ref{wg04});
loosely speaking, it ensures convergence of the (oriented) normal down to (spatial)
scales of $O(\sqrt{h})$. In particular, it rules out ghost interfaces.
Since it will also be used for the variational interpolation,
cf.~Lemma \ref{Le1}, it is formulated for a $[0,1]$-valued sequence $\{u_h\}_{h\downarrow0}$.

\begin{lemma}\label{Le3}
i) Consider a sequence $\{\chi_h\}_{h\downarrow0}$ of $\{0,1\}$-valued functions
on $(0,T)\times[0,1)^d$ that satisfies
\begin{align}\label{wg79}
{\rm esssup}_{t\in(0,T)} E_h(\chi_h(t))
+\int_0^T\frac{1}{2h^2}d_h^2(\chi_h(t),\chi_h(t-h))dt\nonumber\\
\mbox{stays bounded as}\;h\downarrow0,
\end{align}
and that is piecewise constant in the sense of (\ref{wg61}).
Such a sequence is compact in $L^1((0,T)\times[0,1)^d)$; any (weak) limit 
$\chi$ is such that $\nabla\chi$ is a bounded measure, equi-integrable in $t$, with
\begin{align}\label{wg84}
c_0\int|\nabla\chi|dt\le\liminf_{h\downarrow0}\int_0^TE_h(\chi_h(t))dt.
\end{align}
ii) Consider a sequence $\{u_h\}_{h\downarrow 0}$ of $[0,1]$-valued functions
on $(0,T)\times[0,1)^d$ and a $\{0,1\}$-valued function $\chi$ on
$(0,T)\times[0,1)^d$ that satisfies (\ref{wg07}) and (\ref{wg04}) 
(with $\chi_h$ replaced by $u_h$) and
\begin{align}\label{wg94}
{\rm ess\,sup}_{t\in(0,T)} E_h(u_h(t))\quad\mbox{stays bounded as}\;h\downarrow0.
\end{align} 
Then, as measures on $(z,t,x)$-space, we have the weak convergences 
\begin{align}
G_1(z)\frac{1}{\sqrt{h}}u_h(t,x)(1-u_h)(t,x-\sqrt{h}z)dx&dtdz\nonumber\\
\rightharpoonup G_1(z)(\nu\cdot z)_+|\nabla\chi|&dtdz,\label{wg11}\\
G_1(z)\frac{1}{\sqrt{h}}(1-u_h)(t,x)u_h(t,x-\sqrt{h}z)dx&dtdz\nonumber\\
\rightharpoonup G_1(z)(\nu\cdot z)_-|\nabla\chi|&dtdz.\label{wg12}
\end{align}
The test functions may even have polynomial growth in $z$.
\end{lemma}


The next two propositions are at the core and provide the link 
between (\ref{wg67}) and (\ref{wg62}). Proposition \ref{Pr2} ensures that
the metric $d_h$, cf.~(\ref{wg68}), is strong enough to control
the right notion of energy of curves in configuration space.
Proposition \ref{Pr3} makes sure that it is not too strong so that
the metric slope $|\partial E_h|$, cf.~(\ref{wg26}), controls
the gradient of the limiting functional.

\begin{proposition}\label{Pr2}
Suppose that (\ref{wg07}) and the conclusion of 
Lemma \ref{Le3} hold (with $u_h$ replaced by $\chi_h$).
Provided the l.~h.~s.~of (\ref{wg51}) is finite,
there exists $V\in L^2(|\nabla\chi|dt)$
that is the normal velocity in the sense of
\begin{align}\label{wg57}
\partial_t\chi=V|\nabla\chi|\quad\mbox{distributionally},
\end{align}
and that is dominated in the sense of
\begin{align}\label{wg51}
\liminf_{h\downarrow0}\int_0^T\frac{1}{2h^2}d_h^2(\chi_h(t+h),\chi_h(t))dt
\ge c_0\int V^2|\nabla\chi|dt.
\end{align}
\end{proposition}

\begin{proposition}\label{Pr3}
Suppose that the conclusions of Lemma \ref{Le3} ii) hold. 
Then there exists $H\in L^2(|\nabla\chi|dt)$
that is the mean curvature in the sense of (\ref{wg60})
%
%
and that is dominated in the sense of
\begin{align}\label{wg74}
\liminf_{h\downarrow0}\int_0^T\frac{1}{2}|\partial E_h(u_h(t))|^2dt
\ge c_0\int (\frac{H}{2})^2|\nabla\chi|dt.
\end{align}
\end{proposition}


\section{Proofs}

We will repeatedly use the (parabolic) scaling of $G_h$, cf.~(\ref{wg76}),
\begin{align}\label{wg05}
	G_h(z)=\frac{1}{\sqrt{h}^d}G_1(\frac{z}{\sqrt{h}})
\end{align}
and its semi-group property in form of
\begin{align}\label{wg38}
	G_{h}*G_{h'}=G_{h+h'}\quad\mbox{in particular}\quad G_\frac{h}{2}*G_\frac{h}{2}=G_h.
\end{align}
The constant $c_0=\frac{1}{\sqrt{2\pi}}$ appears because of the identity
\begin{align}\label{wg03}
	\int G_1(z)(z_1)_+dz=G^{d=1}_1(0)=c_0,
\end{align}
where $G^{d=1}_1(z_1):=\frac{1}{\sqrt{2\pi}}\exp(-\frac{z_1^2}{2})$ 
denotes the standard Gaussian in a single variable. Indeed, by the factorization
of the $d$-dimensional standard Gaussian into $G_1^{d=1}$ and the
$(d-1)$-dimensional one, and by the normalization of the latter, the
integral in (\ref{wg03}) reduces to $\int_0^\infty G_1^{d=1}z_1dz_1$.
The formula then follows from writing $z_1 G_1^{d=1}$ $=-\frac{d}{dz_1}$ $G_1^{d=1}$.

\medskip

{\sc Proof of Theorem \ref{Th}}.\nopagebreak

Note that Lemma \ref{Le2} allows to make use of Lemma \ref{Le1},
so that we have (\ref{wg67}) with $(E,d,\chi,u)$ replaced by $(E_h,d_h,\chi_h,u_h)$. 
We start with the l.~h.~s.~of (\ref{wg67}), for which we plainly have
\begin{align}\label{wg71}
E_h(\chi^0)\le c_0\int|\nabla\chi^0|.
\end{align}
Indeed, dropping the index $0$, this follows by making the l.~h.~s.~explicit
$\frac{1}{\sqrt{h}}\int G_h(z)(1-\chi)\chi(\cdot-z)dxdz$, which by (\ref{wg05})
and $\chi\in\{0,1\}$ coincides with
$\int G_1(z)\frac{1}{\sqrt{h}}(\chi-\chi(\cdot-\sqrt{h}z))_-dxdz$.
It remains to appeal to the mean-value inequality 
$\int\frac{1}{\sqrt{h}}(\chi-\chi(\cdot-\sqrt{h}z))_-dx$ $\le\int(z\cdot\nu)_-|\nabla\chi|$
and to (\ref{wg03}).

\medskip

Note that because of (\ref{wg67}) and (\ref{wg71}), (\ref{wg79}) is satisfied. 
Hence we may apply Lemma \ref{Le3} i), which yields 
$\chi\in\{0,1\}$ a.~e.~and that $\nabla\chi$
is a bounded measure which is equi-integrable in $t$.
By Lemma \ref{Le3} ii), in view
of the theorem's assumption (\ref{wg04}), we obtain (\ref{wg11}) \& (\ref{wg12})
with $u_h$ replaced by $\chi_h$, so that we may apply Proposition \ref{Pr2}.
We now argue that (\ref{wg69}) \& (\ref{wg70})
imply that (\ref{wg07}) \& (\ref{wg04})
hold with $\chi_h$ replaced by $u_h$, so that we may use Proposition \ref{Pr3} also for $u_h$.
Indeed, (\ref{wg04}) for $u_h$ follows immediately from (\ref{wg04}) for $\chi_h$
and (\ref{wg70}). We now turn to (\ref{wg07}); because it is $[0,1]$-valued,
the sequence $\{u_h\}_{h\downarrow0}$ always admits a subsequence that has a
weak limit $u$, so that it remains to argue that $u=\chi$, w.~l.~o.~g.~assuming
that the entire sequence converges. We momentarily fix $h_0>0$ and note that by 
(\ref{wg38}) together with Jensen's inequality we have for all $h\le 2 h_0$ that
\begin{align*}
\lefteqn{\int|G_{h_0}*(u_h-\chi_h)|^2dxdt\le\int|G_\frac{h}{2}*(u_h-\chi_h)|^2dxdt}\nonumber\\
&\stackrel{(\ref{wg18})}{=}\frac{1}{2\sqrt{h}}\int_0^Td_h^2(u_h(t),\chi_h(t))dt
\stackrel{(\ref{wg69})}{\le}h\sqrt{h}E_h(\chi^0)
\stackrel{(\ref{wg71})}{\le}c_0h\sqrt{h}\int|\nabla\chi^0|,
\end{align*}
so that by lower-semi continuity of the l.~h.~s.~under weak convergence we
obtain $\int|G_{h_0}*(u-\chi)|^2dxdt=0$. From letting $h_0$ tend to zero
we obtain the desired $u=\chi$.

\medskip

Momentarily setting $\rho(t)$ $:=E_h(\chi_h(t))$ $+\frac{1}{2}\int_0^t$ 
$\big(\frac{1}{h^2}d_h^2(\chi_h(s),\chi_h(s-h))$ $+|\partial E_h(u_h(s))|^2\big)ds$,
we note that by definition $\rho(t) =\rho(nh) + \delta(t)$, for
$t\in\big[nh, (n+1) h\big)$, where $\delta(t) := \frac12\int_{nh}^t \big(\frac{1}{h^2}d_h^2(\chi_h(s),\chi_h(s-h))$ $+|\partial E_h(u_h(s))|^2\big)ds $. 
By (\ref{wg67}), $\int_0^T \delta(t)dt$ $ \le h E_h(\chi^0)$.   Hence if we multiply (\ref{wg67}) in form of
$\rho(nh)\le E_h(\chi^0)$ with $\eta(nh)-\eta((n+1)h)$ for some
non-increasing $\eta\in C_0^\infty([0,T))$, we obtain
$\int_0^\infty(-\frac{d\eta}{dt})\rho dt$ $\le \big( \eta(0)+ h \sup \big| \frac{d \eta}{dt}\big| \big) E_h(\chi^0) $.
By an integration by parts and with the choice
$\eta(t) =\max\{\min\{\frac{T-t}{\tau},1\},0\}$, this turns into
%
%
%
%
%
%
\begin{align*}
\lefteqn{\frac{1}{\tau}\int_{T-\tau}^TE_h(\chi_h(t))dt}\nonumber\\
	&+\frac{1}{2}\int_0^{T-\tau} \!\! \big(\tfrac{1}{h^2}d_h^2(\chi_h(t),\chi_h(t-h))
	+|\partial E_h(u_h(t))|^2\big)dt\le \big(1+\tfrac{h}{\tau} \big) E_h(\chi^0).
\end{align*}
Passing to the limit $h\downarrow0$, for the first l.~h.~s.~term, we appeal to (\ref{wg84}) in Lemma \ref{Le3}
with $(0,T)$ replaced by $(T-\tau,T)$.
%
%
For the second l.~h.~s.~term, we apply (\ref{wg51}) (note that we may
extend the integral down to 0 because of the second item in (\ref{wg61}))
and (\ref{wg74}), both with $(0,T)$ replaced by $(0,T-\tau)$. For the r.~h.~s.~term,
we use (\ref{wg71}). Summing up, we obtain
\begin{align*}
	\lefteqn{\frac{c_0}{\tau}\int_{(T-\tau,T)\times[0,1)^d}|\nabla\chi|dt}\nonumber\\
	&+c_0\int_{(0,T-\tau)\times[0,1)^d}\big(V^2+(\frac{H}{2})^2\big)|\nabla\chi|dt
	\le c_0\int|\nabla\chi^0|.
\end{align*}
Dividing by $c_0$ and letting $\tau\downarrow0$ yields (\ref{wg62}).

\medskip

Finally, we argue why (\ref{wg57}) is sufficient to infer (\ref{wg57bis}).
Indeed, by the trivial extension of $\chi_h$ to $t\le 0$, cf.~(\ref{wg61}),
the assumptions (\ref{wg07}) \& (\ref{wg04}) extend to $(-T,T)$, where
for (\ref{wg04}) we appeal to (\ref{wg71}). Likewise,
the l.~h.~s.~integral in (\ref{wg51}) extends to $(-T,T)$. Hence (\ref{wg57})
holds distributionally on $(-T,T)\times[0,1)^d$, which turns into (\ref{wg57bis})
because of $\chi=\chi^0$ for $t<0$.

\bigskip


{\sc Proof of Lemma \ref{Le1}}.

We reproduce the proof of \cite[Theorem 3.1.4 \& Lemma 3.1.3]{AmbrosioGigliSavare}.
We start with the definition of the variational interpolation $u$.
Since by assumption, $({\mathcal M},d)$ is compact and $E$ continuous,
for any $n\in\mathbb{N}$ and any $t\in((n-1)h,nh]$, there exists $u(t)$ that
minimizes
\begin{align}\label{wg96}
	\frac{d^2(u,\chi^{n-1})}{2(t-(n-1)h)}+E(u)\quad\mbox{among}\;u\in{\mathcal M}.
\end{align}
W.~l.~o.~g.~we may assume that $u(nh)=\chi^n$, cf.~(\ref{wg66}), 
so that $u$ is indeed an interpolation
of $\{\chi^n\}_{n\in\mathbb{N}}$. Since by comparison with $u=\chi^{n-1}$
we have $E(u(t))\le E(\chi^{n-1})$, (\ref{wg70}) follows immediately from
the way we defined the piecewise linear interpolation, cf.~(\ref{wg61}).

\medskip

Fixing $n\in\mathbb{N}$ and introducing
\begin{align}\label{wg99}
	e(t):=\min_{u\in{\mathcal M}}\big(\frac{d^2(u,\chi^{n-1})}{2(t-(n-1)h)}+E(u)\big),
	\quad(n-1)h<t\le nh,
\end{align}
we now establish the two crucial inequalities
\begin{align}\label{wg98}
\lefteqn{\frac{d^2(u(s),\chi^{n-1})}{2(s-(n-1)h)(t-(n-1)h)}
	\le\frac{e(s)-e(t)}{t-s}}\nonumber\\
	&\le\frac{d^2(u(t),\chi^{n-1})}{2(s-(n-1)h)(t-(n-1)h)},
\quad(n-1)h<s<t\le nh.
\end{align}
For notational simplicity we consider the case $n=1$; for any $s,t>0$ we have
by definitions (\ref{wg96}) and (\ref{wg99})
\begin{align*}
e(s)&=\frac{1}{2s}d^2(u(s),\chi^0)+E(u(s))\\
&\le \frac{1}{2s}d^2(u(t),\chi^0)+E(u(t))
=(\frac{1}{2s}-\frac{1}{2t})d^2(u(t),\chi^0)+e(t).
\end{align*}
Writing $\frac{1}{2s}-\frac{1}{2t}$ $=\frac{t-s}{2st}$ this gives the upper bound
in (\ref{wg98}) after division by $t-s>0$. 
Exchanging the roles of $s$ and $t$, we likewise get the lower one.

\medskip

We now argue that
\begin{align}\label{wg97}
	|\partial E(u(t))|\le\frac{d(u(t),\chi^{n-1})}{t-(n-1)h}
\quad\mbox{for}\;t\in((n-1)h,nh].
\end{align}
Again, for notational simplicity we consider $n=1$ and give ourselves
a $v\in{\mathcal M}$. By the characterizing property (\ref{wg96})
of $u(t)$ we have $\frac{1}{2t}d^2(u(t),\chi^0)+E(u(t))$
$\le\frac{1}{2t}d^2(v,\chi^0)+E(v)$, so that
$E(u(t))-E(v)$ $\le\frac{1}{2t}(d(v,\chi^0)-d(u(t),\chi^0))(d(v,\chi^0)+d(u(t),\chi^0))$.
By the triangle inequality, this implies
\begin{align*}
	E(u(t))-E(v)\le d(v,u(t))\frac{1}{t}\big(d(u(t),\chi^0)+\frac{1}{2}d(v,u(t))\big),
\end{align*}
so that (\ref{wg97}) follows from definition (\ref{wg26}) of the metric slope.

\medskip

We now may conclude on (\ref{wg67}). By telescoping and 
according to the piecewise constant interpolation,
it is sufficient to establish
\begin{align*}
	E(\chi^n)+\frac{1}{2h}d^2(\chi^n,\chi^{n-1})+\int_{(n-1)h}^{nh}\frac{1}{2}|\partial E(u(s))|^2ds
\le E(\chi^{n-1}),
\end{align*}
which according to (\ref{wg97}) follows from
\begin{align*}
E(\chi^n)+\frac{1}{2h}d^2(\chi^n,\chi^{n-1})
	+\int_{(n-1)h}^{nh}\frac{d^2(u(s),\chi^{n-1})}{2(s-(n-1)h)^2}ds\le E(\chi^{n-1}),
\end{align*}
and with help of (\ref{wg99}) may be rewritten as
\begin{align}\label{so01}
	e(nh)+\int_{(n-1)h}^{nh}\frac{d^2(u(s),\chi^{n-1})}{2(s-(n-1)h)^2}ds
\le E(\chi^{n-1}).
\end{align}
Here comes the argument for (\ref{so01}): We first learn from (\ref{wg98}) that
\begin{align}\label{so02}
((n-1)h,nh]\ni s\mapsto d^2(u(s),\chi^{n-1})\quad\mbox{is monotone increasing}
\end{align}
and thus continuous outside of a countable set of $s$'s. We then learn that
$e$ is locally Lipschitz continuous on $((n-1)h,nh]$ and differentiable
where (\ref{so02}) is continuous. In particular, we have in those (Lebesgue)
almost every time points $s$, $\frac{de}{dt}(s)$ $=-\frac{1}{2s^2}d^2(u(s),\chi^{n-1})$.
Integrating this relationship from some $t\in((n-1)h,nh]$ to $nh$ we obtain
$e(nh)$ $+\int_t^{nh}\frac{1}{2s^2}d^2(u(s),\chi^{n-1})ds$ $\le e(t)$.
Using the obvious $e(t)\le E(\chi^{n-1})$, cf.~(\ref{wg99}), and letting $t\downarrow(n-1)h$
we obtain (\ref{so01}) by monotone convergence.

\medskip

We finally turn to (\ref{wg69}). According to (\ref{so02}) we have
$d^2(u(t),\chi^{n-1})$ $\le d^2(\chi^n,\chi^{n-1})$ for $t\in((n-1)h,nh]$ and thus
$\int_0^\infty d^2(u(t),\chi(t))dt$ $\le\int_0^\infty d^2(\chi(t),\chi(t-h))dt$,
so that (\ref{wg69}) follows from (\ref{wg67}).

\bigskip


{\sc Proof of Lemma \ref{Le2}}.

By the definitions (\ref{wg06}) and (\ref{wg18}), the latter in
conjunction with (\ref{wg38}), we have
\begin{align*}
\frac{1}{2h}d_h^2(u,\chi^{n-1})+E_h(u)=
	\langle u-\chi^{n-1},u-\chi^{n-1}\rangle +\langle 1-u,u\rangle,
\end{align*}
where we momentarily introduced the bilinear form $\langle u,u'\rangle$ 
$:=\frac{1}{\sqrt{h}}\int u$ $ G_h*u'dx$. Since this form is symmetric, we may rewrite
the r.~h.~s.~as $\langle u,1-2\chi^{n-1}\rangle$ $+\langle\chi^{n-1},\chi^{n-1}\rangle$,
so that
\begin{align*}
\frac{1}{2h}d_h^2(u,\chi^{n-1})+E_h(u)=
	\frac{1}{\sqrt{h}}\int u (1-2G_h*\chi^{n-1})+C,
\end{align*}
where $C:=\langle\chi^{n-1},\chi^{n-1}\rangle$ does not depend on $u$. 
It is now obvious that (\ref{wg59}) minimizes $\frac{1}{2h}d_h^2(u,\chi^{n-1})+E_h(u)$
among all $u\in{\mathcal M}$, cf.~(\ref{wg92}).

\medskip

It remains to argue that the metric space $({\mathcal M},d_h)$ is compact
and $E_h$ continuous. Both follows from the fact that $d_h$ metrizes
weak convergence on ${\mathcal M}\subset L^2([0,1)^d)$. The latter can be
seen as follows: In terms of Fourier series, we have
$\frac{2}{\sqrt{h}}d_h^2(u,u')$ 
$=\sum_{k\in2\pi\mathbb{Z}^d}\exp(-\frac{|hk|^2}{2})$ $|{\mathcal F}(u-u')|^2(k)$,
and note $|{\mathcal F}(u-u')|^2(k)$ $\le\int(u-u')^2dx$ $\le 1$.
Hence by dominated convergence,
$d_h(u_n,u)\rightarrow 0$ is equivalent to ${\mathcal F}(u_n-u)(k)\rightarrow 0$
for all $k\in2\pi\mathbb{Z}^d$, which by the $L^2$-boundedness of $\{u_n\}_{n\uparrow\infty}$
$\subset{\mathcal M}$ is equivalent to weak convergence.

\bigskip


{\sc Proof of Lemma \ref{Le3}}. 

\newcounter{Le3S} 
\refstepcounter{Le3S} 

{\sc Step} \arabic{Le3S}.\label{Le3S1}\refstepcounter{Le3S} 
Some useful inequalities on $E_h$ and $d_h$.
We claim for any $[0,1]$-valued function $u$ of space:
\begin{align}
	\int|u-G_h*u|dx
	&\le 2\sqrt{h}E_h(u),\label{wg85}\\
	E_{h_0}(u)&\le E_{h}(u)\quad\mbox{for}\;h_0\in\mathbb{N}^2h,\label{wg86}
\end{align}
which we claim combine to
\begin{align}\label{wg93}
\int(u-G_{h_0}*u)^2dx \le 4\sqrt{h_0}E_{h}(u)\quad\mbox{for all}\;h_0\ge h.
\end{align}
We also claim for any pair of $\{0,1\}$-valued functions $\chi,\chi'$ of space
\begin{align}
\int|\chi-\chi'|dx&\le \frac{1}{2\sqrt{h}}d_h^2(\chi,\chi')
+2\sqrt{h}\big(E_h(\chi)+E_h(\chi')\big).\label{wg83}
\end{align}
We first tackle (\ref{wg85}); by Jensen's inequality in form of
$|u-G_h*u|(x)$ $\le\int G_h(z)|u(x)-u(x-z)|dz$, and by the definition (\ref{wg06})
of $E_h$ which together with the symmetry of $G_h*$ yields
$2\sqrt{h}E_h(u)$ $=\int((1-u)G_h*u$ $+uG_h*(1-u))dx$, (\ref{wg85}) follows 
from the elementary inequality
\begin{align}\label{wg95}
	|u-u'|\le (1-u)u'+u(1-u')\quad\mbox{for}\;u,u'\in[0,1].
\end{align}
\medskip

We now turn to (\ref{wg86}) which we (iteratively) establish in the more general form of
\begin{align}\label{wg91}
	\sqrt{h_0}E_{h_0}\le\sqrt{h}E_{h}
	+\sqrt{h'}E_{h'}\quad\mbox{provided}\;\sqrt{h_0}=\sqrt{h}+\sqrt{h'}.
\end{align}
Indeed, by definition (\ref{wg06}) and the scaling (\ref{wg05})
we have $\sqrt{h}E_h(u)$ $=\int G_1(z)(1-u)(x)u(x-\sqrt{h}z)dxdz$
and by a change of variables in $x$, $\sqrt{h'}E_{h'}(u)$ 
$=\int G_1(z)(1-u)(x-\sqrt{h}z)u(x-(\sqrt{h}+\sqrt{h'})z)dxdz$.
Hence (\ref{wg91}) follows from the elementary inequality
\begin{align*}
	(1-u)u''\le (1-u)u'+(1-u')u''\quad\mbox{for}\;u,u',u''\in[0,1].
\end{align*}
This is equivalent to $u'(u+u''-1)$ $\le uu''$, which because of $u'\in[0,1]$ and $uu''\ge0$
follows from $u+u''-1\le uu''$. The latter is equivalent to $u''(1-u)\le 1-u$,
which holds because of $u,u''\in[0,1]$.
\medskip

We now turn to the upgrade (\ref{wg93}). We first observe that the 
l.~h.~s.~is monotone increasing in $h_0$, as can be seen by the Fourier representation
$\sum_{k}(1-\exp(-\frac{|h_0k|^2}{2}))|{\mathcal F} u(k)|^2$, where ${\mathcal F}u(k)$,
$k\in2\pi\mathbb{Z}^d$, denotes the Fourier series of $u$ and $\exp(-\frac{|h_0k|^2}{2})$
is the Fourier transform of $G_{h_0}$. Now given $h_0\ge h$, 
we write $\sqrt{h_0}=N\sqrt{h}-s$ with $N\in\mathbb{N}$ and $s\in[0,\sqrt{h})$. By the
above monotonicity we have
$\int(u-G_{h_0}*u)^2$ $\le\int(u-G_{N^2h}*u)^2$ $\le\int|u-G_{N^2h}*u|$, 
to which we first apply (\ref{wg85})
with $h$ replaced by $N^2h$ and second apply (\ref{wg86}) with $N^2h$ playing the role of $h_0$.
Hence we end up with $\int(u-G_{h_0}*u)^2$ $\le 2N\sqrt{h}E_h(u)$, which because
of $N\sqrt{h}$ $\le\sqrt{h_0}+\sqrt{h}$ turns into (\ref{wg93}).

\medskip

We finally address (\ref{wg83}), which according to the definitions (\ref{wg06})\&(\ref{wg68})
of $E_h$ and $d_h$ and the simple estimate \eqref{wg85} follows from integrating the following
inequality in $x$, appealing to the symmetry of $G_h*$,
\begin{align*}
|\chi-\chi'|\le(\chi-\chi')G_h*(\chi-\chi')
&+|\chi -G_h*\chi|+|\chi' -G_h*\chi'|.
\end{align*}
%
Writing $|\chi-\chi'| = (\chi-\chi')^2 = (\chi-\chi')G_h*(\chi-\chi') +(\chi-\chi') (\chi-G_h*\chi) +(\chi'-\chi) (\chi'-G_h*\chi')$, we see that the inequality relies on $(\chi-\chi') (\chi-G_h*\chi) \le |\chi-G_h*\chi|$ and on the same inequality with the roles of $\chi$ and $\chi'$ exchanged.

\medskip

{\sc Step} \arabic{Le3S}.\label{Le3S1bis}\refstepcounter{Le3S} 
Modulus of continuity in time: 
%
%
We claim that for every $\{0,1\}$-valued function $\chi$ of time-space
that is piecewise constant in time, cf.~(\ref{wg61}), and any time shift $s\in[0,1]$ we have
\begin{align}\label{wg89}
	I(s)\le C_0\left\{\begin{array}{ccc}
	\frac{s}{\sqrt{h}}&\mbox{for}&s\le h\\
	2\sqrt{h}          &\mbox{for}&h\le s\le\sqrt{h}\\
	4s                 &\mbox{for}&\sqrt{h}\le s
\end{array}\right\}\le 4C_0\sqrt{s},
\end{align}
where in this step of the proof, we use the abbreviation
\begin{align}
I(s)&:=\int_{(s,T)\times[0,1)^d}|\chi(t,x)-\chi(t-s,x)|dxdt,\nonumber\\
C_0&:=\int_h^T\frac{1}{2h^2}d_h^2(\chi(t),\chi(t-h))dt + 4\int_0^TE_h(\chi(t))dt.\label{wg90}
\end{align}
Indeed, for $s\ge h$, we use (\ref{wg83}) with $(\chi,\chi')$ $=(\chi(t),\chi(t-s))$
and integrate in $t\in(s,T)$ to obtain
\begin{align}\label{wg87}
I(s)\le\int_s^T\frac{1}{2\sqrt{h}}d_h^2(\chi,\chi(\cdot-s))dt
+4\sqrt{h}\int_0^TE_h(\chi)dt,
\end{align}
where we write $\chi(\cdot-s)$ for the time-shifted function $(t,x)\mapsto\chi(t-s,x)$.

We first use this to treat (\ref{wg89}) in case of $s\le h$: By the piecewise constant interpolation,
cf.~(\ref{wg61}), we have $I(s)\le \frac{s}{h}I(h)$,
into which we insert (\ref{wg87}) for $s=h$ in form of $I(h)\le C_0\sqrt{h}$.
We now treat (\ref{wg89}) in case of $h\le s\le\sqrt{h}$, and first restrict
ourselves to $s=Nh$ with $N\in\mathbb{N}$
in order to use the triangle inequality for $d_h$ in form of
$d_h^2(\chi,\chi(\cdot-s))$ $\le N\sum_{n=1}^Nd_h^2(\chi(\cdot-(n-1)h),\chi(\cdot-nh))$
so that we obtain from (\ref{wg87})
\begin{align}\label{wg88}
I(s)&\le(\frac{s}{h})^2\int_h^T\frac{1}{2\sqrt{h}}d_h^2(\chi,\chi(\cdot-h))dt
+4\sqrt{h}\int_0^TE_h(\chi)dt\nonumber\\
&\stackrel{(\ref{wg90})}{\le} C_0\max\{\frac{s^2}{\sqrt{h}},\sqrt{h}\}= C_0\sqrt{h}.
\end{align}
For the unrestricted range of $h\le s\le\sqrt{h}$, we write $s=Nh+s'$ with $N\in\mathbb{N}$ and
$s'\in[0,h)$, use $I(s)\le I(Nh)+I(s')$, and
appeal to (\ref{wg88}) for the first contribution and to (\ref{wg89}) in the
previously treated case of $s'\le h$ for the second contribution.
Finally, for (\ref{wg89}) in the remaining case of $s\ge\sqrt{h}$, we write $s=N\sqrt{h}+s'$ 
with $N\in\mathbb{N}$ and $s'\in[0,\sqrt{h})$, use $I(s)\le N I(\sqrt{h})+I(s')$, and
appeal to the previously treated case of (\ref{wg89}) for both terms.

\medskip

{\sc Step} \arabic{Le3S}.\label{Le3S3}\refstepcounter{Le3S}
Proof of the compactness statement. From (\ref{wg93}) and thanks
to the first part of our assumed bound (\ref{wg79}), we learn that
$G_{h_0}*\chi_h$ is close to $\chi_h$ in $L^\infty((0,T), L^2([0,1)^d))$
(and thus in $L^1((0,T)\times[0,1)^d)$), as $h_0\downarrow 0$,
uniformly in $h\downarrow0$. Hence it remains to argue for fixed $h_0>0$
that $\{G_{h_0}*\chi_h\}_{h\downarrow0}$ is compact in $L^1$. Because of
the convolution in space, and of the equi-integrability following 
from $G_{h_0}*\chi_h$ $\in [0,1]$,
this follows from a modulus of continuity in time in $L^1$
that is uniform in $h\downarrow0$. Thanks to our assumption (\ref{wg79}),
this holds for $\chi_h$ itself by Step \ref{Le3S1bis}; it transmits to $G_{h_0}*\chi_h$
by Jensen's inequality.

\medskip

{\sc Step} \arabic{Le3S}.\label{Le3S3bis}\refstepcounter{Le3S}
Before establishing the exact inequality (\ref{wg84}), which will be done 
at the end of Step \ref{Le3S5},
it is convenient to first argue that $\nabla\chi$ is a bounded measure,
equi-integrable in $t$, under the mere assumption (\ref{wg94}).
We focus on $\partial_1\chi$, give
ourselves a $\zeta\in C^\infty_0((0,T)\times[0,1)^d)$ and note that
as a consequence of (\ref{wg03}) we have
\begin{align*}
-c_0\partial_1\zeta
&=\lim_{h\downarrow0}\frac{1}{\sqrt{h}}\int_{z_1>0}(\zeta-\zeta(\cdot+\sqrt{h}z))G_1(z)dz\nonumber\\
&\stackrel{(\ref{wg05})}{=}\lim_{h\downarrow0}
\frac{1}{\sqrt{h}}\int_{z_1>0}(\zeta-\zeta(\cdot+z))G_h(z)dz
\end{align*}
uniformly in time-space, where we write $\zeta(\cdot+z)$ for the space-shifted function
$(t,x)\mapsto\zeta(t,x+z)$.
Hence it follows from (\ref{wg07}) that
\begin{align*}
-c_0\int\partial_1\zeta\chi dxdt
	=\lim_{h\downarrow0}\int_{z_1>0}\zeta \frac{1}{\sqrt{h}}(\chi_h-\chi_h(\cdot-z))G_h(z) dxdtdz
\end{align*}
and thus
\begin{align*}
c_0\big|\int\partial_1\zeta\chi dxdt\big|
	\le\int_0^T\sup_x|\zeta|dt\liminf_{h\downarrow0}
	\mathop{\rm ess\,sup}_t\frac{1}{\sqrt{h}}\int|\chi_h-G_h*\chi_h|dx.
\end{align*}
According to (\ref{wg85}), the second r.~h.~s.~factor is estimated by
$2\mathop{\rm ess\,sup}_t$ $ E_h(\chi_h)$, the $\liminf_{h\downarrow0}$ of which is bounded
by our assumption (\ref{wg94}).

\medskip

{\sc Step} \arabic{Le3S}.\label{Le3S4}\refstepcounter{Le3S}
Turning to (\ref{wg11}) \& (\ref{wg12}), it is convenient to have
the following equi-integrability of the non-negative 
\begin{align*}
	\rho_h(z,t,x):=G_1(z)\frac{1}{\sqrt{h}}\big(&(1-u_h)(t,x)u_h(t,x-\sqrt{h}z)\nonumber\\
	&+u_h(t,x)(1-u_h)(t,x-\sqrt{h}z)\big)
\end{align*}
in the (non-compact) variables $t$ and $z$ in the sense of
\begin{align*}
	\int\exp(\frac{3|z|^2}{8})\rho_h(z,t,x)dxdz\le 2^{d+2}E_h(u_h(t)),
\end{align*}
cf.~(\ref{wg94}). Indeed, we first observe that $G_4(z)$ 
$=\frac{1}{2^d}\exp(\frac{3|z|^2}{8})G_1(z)$, so that by the scaling (\ref{wg05})
we obtain $\int\exp(\frac{3|z|^2}{8})\rho_h dxdz$ $=\frac{2^d}{\sqrt{h}}$ 
$\int\big((1-u_h)$ $G_{4h}*u_h$ $+u_hG_{4h}*(1-u_h)\big)dx$, so that by symmetry of $G_{4h}$
and the definition (\ref{wg06}) of $E_{4h}$,
$\int\exp(\frac{3|z|^2}{4})\rho_h dxdz$ $=2^{d+2}E_{4h}(u_h)$,
so that it remains to appeal to (\ref{wg86}).

\medskip

{\sc Step} \arabic{Le3S}.\label{Le3S5}\refstepcounter{Le3S}
Proof of (\ref{wg11}) \& (\ref{wg12}).
We focus on (\ref{wg11}); according to Step \ref{Le3S4}, it is sufficient
to treat bounded continuous test functions $\zeta(z,t,x)$, which by linearity
and splitting into positive and negative part we may assume to be $[0,1]$-valued. 
The statement splits into the local lower bound
\begin{align}\label{wg01}
\liminf_{h\downarrow 0}&
\int\zeta G_1(z) \frac{1}{\sqrt{h}}u_h(1-u_h)(\cdot-\sqrt{h}z)dxdtdz\nonumber\\
&\ge\int\zeta G_1(z)(\nu\cdot z)_+|\nabla\chi|dtdz
\end{align}
and the global upper bound
\begin{align}\label{wg02}
\limsup_{h\downarrow 0}&
\int G_1(z)\frac{1}{\sqrt{h}}u_h(1-u_h)(\cdot-\sqrt{h}z)dxdtdz\nonumber\\
&\le\int G_1(z)(\nu\cdot z)_+|\nabla\chi|dtdz.
\end{align}
Indeed, splitting a given test function $\zeta=1-(1-\zeta)$, 
appealing to linearity and using (\ref{wg01}) 
with $\zeta$ replaced by $1-\zeta\in[0,1]$, we obtain also
the local upper bound. 

\medskip

We note that the global upper bound (\ref{wg02}) is nothing else than
our assumption (\ref{wg04}) (with $\chi_h$ replaced by $u_h$): 
For the l.~h.~s.~this follows from the 
the scaling (\ref{wg05}), the symmetry of $G_h*$, and
the definition (\ref{wg06}) of $E_h$. For the r.~h.~s.~this follows from (\ref{wg03}).

\medskip

Hence it remains to establish (\ref{wg01}) by a typical l.~s.~c.~argument:
By Fatou's Lemma, it is enough to establish for fixed $z$
\begin{align*}
\liminf_{h\downarrow 0}
\int\zeta\frac{1}{\sqrt{h}}u_h(1-u_h)(\cdot-\sqrt{h}z)dxdt
\ge\int\zeta(\nu\cdot z)_+|\nabla\chi|dt.
\end{align*}
Since by Step \ref{Le3S3bis}, we already know that $\nabla\chi$ is a bounded
measure, we have $(\nu\cdot z)_+|\nabla\chi|=(\partial_z\chi)_+$. Hence by the definition
of the positive part of a measure, and by the equi-integrability in $t$ established
in Step \ref{Le3S4}, it is enough to establish for any 
non-negative $\zeta\in C_0^\infty((0,T)\times[0,1)^d)$
\begin{align*}
\liminf_{h\downarrow 0}
\int\zeta\frac{1}{\sqrt{h}}u_h(1-u_h)(\cdot-\sqrt{h}z)dxdt
\ge-\int\partial_z\zeta \chi dxdt.
\end{align*}
Since by assumption (\ref{wg07}), the r.~h.~s.~is the limit of
\begin{align*}
\int\frac{1}{\sqrt{h}}(\zeta-\zeta(\cdot+\sqrt{h}z)) u_h
=\int\zeta\frac{1}{\sqrt{h}}(u_h-u_h(\cdot-\sqrt{h}z)),
\end{align*}
the statement follows from the elementary inequality
$u-u'$ $\le u(1-u')$ that is valid for any $u,u'\in[0,1]$.

\medskip

We note that this argument for (\ref{wg01}) did not involve the extra assumption
(\ref{wg04}) and thus applies also under the assumptions of part i)
of the lemma. Statement (\ref{wg01}) applied to $\zeta=1$ yields (\ref{wg84})
(for the same reason, given above, that (\ref{wg02}) is a mere reformulation of (\ref{wg04})).

\bigskip


{\sc Proof of Proposition \ref{Pr2}}.\nopagebreak

Note that by definition (\ref{wg18}) of $d_h$ and (\ref{wg38}) we have
\begin{align*}
\lefteqn{\int_h^T\frac{1}{2h^2}d_h^2(\chi_h(t),\chi_h(t-h))dt}\nonumber\\
&=\int_{(h,T)\times[0,1)^d}\frac{1}{h\sqrt{h}}|G_\frac{h}{2}*(\chi_h-\chi_h(\cdot-h))|^2dxdt,
\end{align*}
which motivates to introduce the non-negative (bounded) ``dissipation measure'' 
on $(0,T)\times[0,1)^d$ (after possibly passing to a subsequence) 
\begin{align}\label{wg35}
\mu:=\lim_{h\downarrow 0}\frac{1}{h\sqrt{h}}|G_\frac{h}{2}*(\chi_h-\chi_h(\cdot-h))|^2.
\end{align}
In fact, we shall establish (\ref{wg51}) in the localized form of
\begin{align}\label{wg53}
c_0V^2|\nabla\chi|dt\le\mu.
\end{align}

\medskip

We now give an outline of the proof, which amounts to a better quantification
of Step \ref{Le3S1bis} in the proof of Lemma \ref{Le3}.
In deriving this estimate on the distribution $\partial_t\chi$,
we work on an intermediate time scale, which we fix to be an (eventually small) fraction
of the characteristic spatial scale:
\begin{align}\label{wg44}
\tau:=\alpha\sqrt{h}\quad\mbox{for some fixed}\;\alpha\in(0,\infty).
\end{align}
We consider the corresponding increments and their positive and negative parts 
\begin{align}\label{wg40}
	\delta\chi:=\chi_h-\chi_h(\cdot-\tau)\quad\mbox{and}\quad
	\delta\chi_\pm:=\max\{0,\pm\delta\chi\}.
\end{align}
%
Using $|\delta\chi| = \delta\chi_+ + \delta\chi_-$ we then write
\begin{align}\label{wg54}
	\frac{1}{\sqrt{h}}|\delta\chi|=\frac{1}{\sqrt{h}}\big(&
 \delta\chi_+ G_h*\delta\chi_+ +\delta\chi_+ G_h*(1-\delta\chi_+)\nonumber\\
	+&\delta\chi_- G_h*\delta\chi_- +\delta\chi_- G_h*(1-\delta\chi_-)\big).
\end{align}
In Step \ref{Pr1S1} we will show, in the sense of closeness between distributions,
\begin{align*}
\frac{1}{\sqrt{h}}|\delta\chi|\approx\frac{1}{\sqrt{h}}\big(
	\delta\chi G_h*\delta\chi +\delta\chi_+ G_h*(1-\delta\chi_+)
	+\delta\chi_- G_h*(1-\delta\chi_-)\big).
\end{align*}
The main idea is to estimate the first r.~h.~s.~by the dissipation measure $\mu$
(in Step \ref{Pr1S2}) and to estimate the two last contributions by
the surface measure $|\nabla\chi|dt$ (in Step \ref{Pr1S3}).
However, this just suffices to show that $\partial_t\chi$
is absolutely continuous w.~r.~t. $|\nabla\chi|dt$, as is carried out at
the beginning of Step \ref{Pr1S4} on the level of $O(\alpha)$ as $\alpha\downarrow0$.
In order to retrieve (\ref{wg53}), we need the finer estimate in Step \ref{Pr1S3}
and look at level $O(\alpha^2)$.

\medskip

\newcounter{Pr1S} 
\refstepcounter{Pr1S} 

{\sc Step} \arabic{Pr1S}.\label{Pr1S1}\refstepcounter{Pr1S}
We claim that the mixed term vanishes distributionally:
\begin{align}\label{wg33}
\lim_{h\downarrow0}\frac{1}{\sqrt{h}}(\delta\chi_+G_h*\delta\chi_-+\delta\chi_-G_h*\delta\chi_+)
=0.
\end{align}
Spelling out the $z$-integral, we want to show that the distributional limit of
\begin{align}\label{wg34}
\int G_h(z)\frac{1}{\sqrt{h}}\big(\delta\chi_+\delta\chi_-(\cdot-z)
+\delta\chi_-\delta\chi_+(\cdot-z)\big)dz
\end{align}
vanishes. Fixing some unit vector $\nu_0$, we split the expression into
\begin{align}\label{wg31}
\int_{\nu_0\cdot z\ge 0} G_h(z)\frac{1}{\sqrt{h}}\big(\delta\chi_+\delta\chi_-(\cdot-z)
	        +\delta\chi_-\delta\chi_+(\cdot-z)\big)dz
\end{align}
and the analogous expression on $\{\nu_0\cdot z\le 0\}$.
We note that by definition of (\ref{wg40}),
we have $\delta\chi_+=\chi_h(1-\chi_h)(\cdot-\tau)$ and
$\delta\chi_-=\chi_h(\cdot-\tau)(1-\chi_h)$ and thus
$\delta\chi_+\delta\chi_-(\cdot-z)$
$=\chi_h(1-\chi_h)(\cdot-\tau)\chi_h(\cdot-\tau-z)(1-\chi_h)(\cdot-z)$
$\le(1-\chi_h)(\cdot-\tau)\chi_h(\cdot-\tau-z)$ and likewise
$\delta\chi_-\delta\chi_+(\cdot-z) =\chi_h(\cdot-\tau) $ $(1-\chi_h)\chi_h(\cdot-z)(1-\chi_h)(\cdot-\tau-z)$
$\le(1-\chi_h)\chi_h(\cdot-z)$. Here, with a slight abuse of notation, $\chi(\cdot-\tau-z) (t,x) = \chi(t-\tau, x-z)$.  Hence the distributional limit of (\ref{wg31})
is dominated by the one of
\begin{align*}
\int_{\nu_0\cdot z\ge 0} G_h(z)\frac{1}{\sqrt{h}}\big(
&(1-\chi_h)(\cdot-\tau)\chi_h(\cdot-\tau-z)\nonumber\\
+&(1-\chi_h)\chi_h(\cdot-z)\big)dz.
\end{align*}
Since the first contribution differs from the second one just by a (vanishing) time shift,
which we may put into the continuous test function,
the distributional limit is equal to the weak limit of
\begin{align*}
2\int_{\nu_0\cdot z\ge 0} G_h(z)\frac{1}{\sqrt{h}}(1-\chi_h)\chi_h(\cdot-z)dz,
\end{align*}
provided the latter exists.
Appealing to the scaling (\ref{wg05}), 
according to (\ref{wg12}) in Lemma \ref{Le3} which we test with the
characteristic function of the closed set $\{\nu_0\cdot z\ge 0\}$, 
the weak limit of this term is dominated by the measure
$2\int_{\nu_0\cdot z\ge0}G_1(z)(\nu\cdot z)_-|\nabla\chi|dtdz$.
Treating the contribution of $\{\nu_0\cdot z\le0\}$ in a similar way (exchanging the roles
of $\chi$ and $1-\chi$), the weak limit of that contribution is dominated by
$2\int_{\nu_0\cdot z\le0}G_1(z)(\nu\cdot z)_+|\nabla\chi|dtdz$.

\medskip

Hence we have shown that the weak limit $\lambda\ge 0$ of  (\ref{wg34})
satisfies 
\begin{align}\label{wg63}
\lambda\le 2\Big(\int_{\nu_0\cdot z\ge0}G_1(\nu\cdot z)_-dz
+\int_{\nu_0\cdot z\le0}G_1(\nu\cdot z)_+dz\Big)|\nabla\chi|dt.
\end{align}
In particular, we have $\lambda\le 4c_0|\nabla\chi|dt$, cf.~(\ref{wg03}), so that 
there is a $\theta\in L^1(|\nabla\chi|dt)$ with $\lambda=\theta|\nabla\chi|dt$,
which allows us the rewrite (\ref{wg63}) as
\begin{align*}
\theta\le2\int_{\nu_0\cdot z\ge0}G_1(\nu\cdot z)_-dz
	+2\int_{\nu_0\cdot z\le0}G_1(\nu\cdot z)_+dz\\
|\nabla\chi|dt-\mbox{a.~e.}\;\mbox{and for all}\;\nu_0\in\mathbb{R}^d.
\end{align*}
An elementary separability argument allows to exchange the order between the $\forall$,
so that we may choose $\nu_0=\nu$, obtaining $\theta\le0$ $|\nabla\chi|dt$-a.~e.~
and thus $\lambda\le 0$, yielding (\ref{wg33}).

\medskip

{\sc Step} \arabic{Pr1S}.\label{Pr1S2}\refstepcounter{Pr1S}
We claim that in a distributional sense
\begin{align}\label{wg37}
\limsup_{h\downarrow 0}\frac{1}{\sqrt{h}}\delta\chi G_h*\delta\chi\le\alpha^2\mu.
\end{align}
Indeed, by definition (\ref{wg35}), we may split this into
\begin{align}
\lim_{h\downarrow 0}\frac{1}{\sqrt{h}}
\big(\delta\chi G_h*\delta\chi-|G_\frac{h}{2}*\delta\chi|^2\big)=0\quad\mbox{and}\label{wg36bis}\\
\limsup_{h\downarrow 0}\Big(\frac{1}{\sqrt{h}}|G_\frac{h}{2}*\delta\chi|^2
	-\alpha^2\frac{1}{h\sqrt{h}}|G_\frac{h}{2}*(\chi_h-\chi_h(\cdot-h))|^2\Big)
	\le 0.\label{wg36}
\end{align}
We start with (\ref{wg36}) and assume w.~l.~o.~g.~that $\tau=Nh$ for some $N\in\mathbb{N}$
so that by telescoping
%
%
and the Cauchy-Schwarz inequality in $n$
\begin{align*}
\lefteqn{\frac{1}{\sqrt{h}}|G_\frac{h}{2}*(\chi_h-\chi_h(\cdot-\tau))|^2}\nonumber\\
&\le N\sum_{n=0}^{N-1}\frac{1}{\sqrt{h}}|G_\frac{h}{2}*(\chi_h(\cdot-nh)-\chi_h(\cdot-(n+1)h))|^2.
\end{align*}
Appealing to $N=\frac{\alpha}{\sqrt{h}}$, the r.~h.~s.~may be rewritten as
\begin{align*}
\alpha^2\frac{1}{N}\sum_{n=0}^{N-1}
\frac{1}{h\sqrt{h}}|G_\frac{h}{2}*(\chi_h(\cdot-nh)-\chi_h(\cdot-(n+1)h))|^2.
\end{align*}
Note that this is an average of time shifts of $\alpha^2
\frac{1}{h\sqrt{h}}|G_\frac{h}{2}*(\chi_h-\chi_h(\cdot-h))|^2$;
because of $Nh=O(\sqrt{h})=o(1)$ all these time shifts are small,
so that the (non-negative) expression has the same (bounded) weak limit as
$\alpha^2\frac{1}{h\sqrt{h}}|G_\frac{h}{2}*(\chi_h-\chi_h(\cdot-h))|^2$ itself.
This yields (\ref{wg36}).

\medskip

We now turn to (\ref{wg36bis}). By the semi-group property (\ref{wg38})
and the symmetry of $G_\frac{h}{2}*$,
we have for a smooth test function $\zeta$
\begin{align*}
\int\zeta\big(\delta\chi G_h*\delta\chi-|G_\frac{h}{2}*\delta\chi|^2\big)
=-\int[\zeta,G_\frac{h}{2}*]\delta\chi\,G_\frac{h}{2}*\delta\chi,
\end{align*}
where $[\zeta,G_\frac{h}{2}*]$ denotes the commutator of multiplying
with $\zeta$ and convolving with $G_\frac{h}{2}$.
Hence by the boundedness of $\frac{1}{\sqrt{h}}|G_\frac{h}{2}*\delta\chi|^2$
as a sequence of measures, which follows from (\ref{wg36}), it is enough 
to establish
\begin{align*}
\lim_{h\downarrow 0}\frac{1}{\sqrt{h}}\int|[\zeta,G_\frac{h}{2}*]\delta\chi|^2=0.
\end{align*}
We spell out the integrand:
\begin{align*}\big([\zeta,G_\frac{h}{2}*]\delta\chi\big)(t,x)
=\int G_\frac{h}{2}(z)(\zeta(t,x)-\zeta(t,x-z))\delta\chi(t,x-z)dz,
\end{align*}
so that $|[\zeta,G_\frac{h}{2}*]\delta\chi|(t,x)$
$\le\sup|\nabla\zeta|\int G_\frac{h}{2}(z)|z||\delta\chi(t,x-z)|dz$
and thus

\begin{align*}
\frac{1}{\sqrt{h}}\int|[\zeta,G_{\frac{h}{2}}*]\delta\chi|^2dx
\le\frac{1}{\sqrt{h}}\big(\sup|\nabla\zeta|\int G_\frac{h}{2}(z)|z|dz\big)^2\int|\delta\chi|^2dxdt.
\end{align*}
By the scaling (\ref{wg05}), the r.~h.~s.~is $O(\sqrt{h})$ and thus vanishing.

\medskip

{\sc Step} \arabic{Pr1S}.\label{Pr1S3}\refstepcounter{Pr1S}
For given unit vector $\nu_0$ and $V_0\in(0,\infty)$ we claim that in a distributional sense
\begin{align}\label{wg55}
\lefteqn{\limsup_{h\downarrow0}\frac{1}{\sqrt{h}}\Big(
\delta\chi_+\,G_h*(1-\delta\chi_+)+\delta\chi_-\,G_h*(1-\delta\chi_-)}\nonumber\\
&-2\int_{z\cdot\nu_0>\alpha V_0}G_1(z)dz|\delta\chi|\Big)
\le2\int_{0\le z\cdot\nu_0\le \alpha V_0}G_1(z)|z\cdot\nu||\nabla\chi|dt.
\end{align}
We split this into three steps: 1) The l.~h.~s.~may be substituted according to
\begin{align}\label{wg47}
\lim_{h\downarrow 0}
\frac{1}{\sqrt{h}}\Big(&
\delta\chi_\pm\,G_h*(1-\delta\chi_\pm)\nonumber\\
&-\int_{z\cdot\nu_0\ge0}G_h(z)|\delta\chi_\pm-\delta\chi_\pm(\cdot-z)|dz\Big)=0.
\end{align}
2) The integrand, which is a second-order difference, satisfies the two inequalities
\begin{align}\label{wg39}
\lefteqn{|\delta\chi_+-\delta\chi_+(\cdot-z)|+|\delta\chi_--\delta\chi_-(\cdot-z)|}\nonumber\\
&\le\left\{\begin{array}{c}
|\chi_h-\chi_h(\cdot-z)|+|\chi_h(\cdot-\tau)-\chi_h(\cdot-\tau-z)|\\
	|\delta\chi|+|\delta\chi(\cdot-z)|
\end{array}\right\},
\end{align}
where the first inequality is ``space-like'' and we use it on the set
$\{0\le z\cdot\nu_0\le\tau V_0\}$, while the second one is ``time-like''
and we use it on the complement $\{z\cdot\nu_0>\tau V_0\}$.
3) We finally argue that 
\begin{align}
	\limsup_{h\downarrow0}&\frac{1}{\sqrt{h}}\int_{0\le z\cdot\nu_0\le\tau V_0}G_h(z)
	\big(|\chi_h-\chi_h(\cdot-z)|+|\chi_h(\cdot-\tau)\nonumber\\
	&-\chi_h(\cdot-z-\tau)|\big)dz
	\le2\int_{0\le z\cdot\nu_0\le\alpha V_0}G_1(z)|z\cdot\nu||\nabla\chi|dtdz,\label{wg43}\\
	\lim_{h\downarrow0}&\frac{1}{\sqrt{h}}\Big(\int_{z\cdot\nu_0>\tau V}G_h(z)
\big(|\delta\chi|+|\delta\chi(\cdot-z)|\big)dz\nonumber\\
&-2\frac{1}{\sqrt{h}}\int_{z\cdot\nu_0>\alpha V} G_1(z)dz
	|\delta\chi|\Big)=0.\label{wg42}
\end{align}

\medskip

We start with (\ref{wg39}); the second inequality is obvious by the
triangle inequality in form of $|\delta\chi_\pm-\delta\chi_\pm(\cdot-z)|$
$\le\delta\chi_\pm+\delta\chi_\pm(\cdot-z)$ and by $\delta\chi_++\delta\chi_-=|\delta\chi|$.
In view of the definition (\ref{wg40}), the first inequality in (\ref{wg39})
follows from the elementary inequality
\begin{align}\label{wg41}
	|(a-a')_+-(b-b')_+|&+|(a-a')_--(b-b')_-|\nonumber\\&\le|a-b|+|a'-b'|,
\end{align}
which can be seen by distinguishing two cases: Case 1): $(a-a')(b-b')\ge 0$,
w.~l.~o.~g. by symmetry $(a,a',b,b')$ $\leadsto(-a,-a',-b,-b')$ we may assume $a-a',b-b'\ge 0$, 
in which case (\ref{wg41}) turns into the obvious $|(a-a')-(b-b')|$ $\le|a-b|+|a'-b'|$.
Case 2): $(a-a')(b-b')\le 0$, by the same symmetry 
we may assume $a-a'\ge0\ge b-b'$, in which 
case (\ref{wg41}) turns into the obvious $(a-a')+(b'-b)\le|a-b|+|a'-b'|$.

\medskip

We now argue for (\ref{wg43}) \& (\ref{wg42}). Putting the vanishing shift $\tau$ in the time
variable on the continuous test function, (\ref{wg43}) follows once we argue that
\begin{align}\label{wg46}
\limsup_{h\downarrow0}&\frac{1}{\sqrt{h}}\int_{0\le z\cdot\nu_0\le\tau V_0}G_h(z)
|\chi_h-\chi_h(\cdot-z)|dz\nonumber\\
&\le\int_{0\le z\cdot\nu_0\le\alpha V_0}G_1(z)|z\cdot\nu||\nabla\chi|dtdz,
\end{align}
which because of non-negativity implies that the l.~h.~s.~admits a limit
as a measure on $(0,T)\times [0,1)^d$. 
Appealing to scaling (\ref{wg05})
in form of $G_h(z)dz=G_1(\frac{z}{\sqrt{h}})d\frac{z}{\sqrt{h}}$ and noting
that $z\cdot\nu\le\tau V_0$ is equivalent to $\frac{z}{\sqrt{h}}\cdot\nu\le\alpha V_0$,
cf.~(\ref{wg44}), (\ref{wg46}) follows from taking the sum of (\ref{wg11}) and (\ref{wg12}),
appealing to (\ref{wg95}),
testing with the characteristic function of the closed set $\{0\le z\cdot\nu_0\le\alpha V_0\}$,
and integrating out $z$.

\medskip

We now turn to (\ref{wg42}) and note that it follows from
\begin{align*}
\lim_{h\downarrow0}\frac{1}{\sqrt{h}}
\int_{z\cdot\nu_0>\tau V}G_h(z)(|\delta\chi(\cdot-z)|-|\delta\chi|)dz=0,
\end{align*}
which testing with $\zeta\in C_0^\infty((0,T)\times[0,1)^d)$ assumes the form 
\begin{align*}
\lim_{h\downarrow0}\frac{1}{\sqrt{h}}
\int_{z\cdot\nu_0>\tau V}G_h(z)(\zeta(\cdot+z)-\zeta)|\delta\chi|dxdtdz=0.
\end{align*}
This holds, since the integral can be estimated by
\begin{align}\label{wg49}
\lefteqn{\sup|\nabla\zeta|\int\frac{|z|}{\sqrt{h}}G_h(z)|\delta\chi|dxdtdz}\nonumber\\
&\stackrel{(\ref{wg05}),(\ref{wg40})}{=}
\sup|\nabla\zeta|\int|z|G_1(z)dz\int|\chi_h-\chi_h(\cdot-\tau)|dxdt,
\end{align}
and the last contribution vanishes in the limit since $\tau$ does, and since 
$\{\chi_h\}_{h\downarrow0}$ is compact
in $L^1((0,T)\times[0,1)^d)$, cf.~part i) of Lemma \ref{Le3}.

\medskip

We finally address (\ref{wg47}) and focus on the $+$-part. We first argue that
\begin{align}\label{wg50}
\lim_{h\downarrow 0}\frac{1}{\sqrt{h}}\big(\delta\chi_+ G_h*(1-\delta\chi_+)
-(1-\delta\chi_+)G_h*\delta\chi_+\big)=0.
\end{align}
Spelling out the $z$-integral, and using that $G_h$ is even, this follows from
\begin{align}\label{wg64}
\lefteqn{\lim_{h\downarrow 0}\frac{1}{\sqrt{h}}
\int G_h(z)\delta\chi_+(1-\delta\chi_+)(\cdot-z)dz}\nonumber\\
&=\lim_{h\downarrow 0}\frac{1}{\sqrt{h}}\int G_h(z)\delta\chi_+(\cdot+z)(1-\delta\chi_+)dz.
\end{align}
Since the second function differs from the first just by a spatial shift of $z$,
the limits coincide provided the l.~h.~s.~limit is finite, which by the non-negativity
of the functions follows if their integral remains bounded.
The l.~h.~s.~integral indeed remains bounded since $\delta\chi_+(1-\delta\chi_+)(\cdot-z)$
$\le\chi_h(1-\chi_h)(\cdot-z)+(1-\chi_h)(\cdot-\tau)\chi_h(\cdot-\tau-z)$, which follows from $\delta \chi_+ = \chi_h(1-\chi_h)(\cdot-\tau)$,
and since the integral of the first summand remains bounded by (\ref{wg11}),
whereas the integral of the second summand is oblivious to the (vanishing) time shift
and thus remains bounded by (\ref{wg12}).
%
%
%
%
%
%

\medskip

Equipped with (\ref{wg50}),
we now may substitute $\frac{1}{\sqrt{h}}\delta\chi_+ G_h*(1-\delta\chi_+)$
by $\frac{1}{2}\frac{1}{\sqrt{h}}\big(\delta\chi_+ G_h*(1-\delta\chi_+)
+(1-\delta\chi_+) G_h*\delta\chi_+\big)$, which we may write as
\begin{align*}
	\frac{1}{2\sqrt{h}}\int G_h(z)|\delta\chi_+-\delta\chi_+(\cdot-z)|dz,
\end{align*}
where we used for any $a,b\in\{0,1\}$ that $a(1-b)+(1-a)b=|a-b|$. Hence in
order to obtain (\ref{wg47}), we need the two function sequences
\begin{align*}
\frac{1}{\sqrt{h}}\int_{\pm z\cdot\nu_0\ge 0} G_h(z)|\delta\chi_+-\delta\chi_+(\cdot-z)|dz
\end{align*}
to have the same limit. Again by the evenness of $G_h$, these two functions only
differ by a spatial shift $z$. The same argument as for (\ref{wg64})
shows that the limits agree.
%
%

\medskip

{\sc Step} \arabic{Pr1S}.\label{Pr1S4}\refstepcounter{Pr1S}
Conclusion. We start from the identity (\ref{wg54})
in form of 
\begin{align*}
\frac{1}{\sqrt{h}}|\delta\chi|&=\frac{1}{\sqrt{h}}\big(
\delta\chi G_h*\delta\chi\nonumber\\
&+\delta\chi_+G_h*\delta\chi_-+\delta\chi_-G_h*\delta\chi_+\nonumber\\
&+\delta\chi_+ G_h*(1-\delta\chi_+)+\delta\chi_- G_h*(1-\delta\chi_-)\big),
\end{align*}
or rather, using $2\int_{z\cdot\nu_0\ge 0}G_1dz=1$,
\begin{align*}
\lefteqn{2\int_{0\le z\cdot\nu_0\le \alpha V_0}G_1dz\frac{1}{\sqrt{h}}|\delta\chi|
=\frac{1}{\sqrt{h}}\big(
	\delta\chi G_h*\delta\chi}\nonumber\\
&+\delta\chi_+G_h*\delta\chi_-+\delta\chi_-G_h*\delta\chi_+\nonumber\\
&+\delta\chi_+ G_h*(1-\delta\chi_+)+\delta\chi_- G_h*(1-\delta\chi_-)
-2\int_{z\cdot\nu_0>\alpha V_0}G_1dz|\delta\chi|\big).
\end{align*}
We note that by an elementary lower-semi-continuity argument
based on the definitions (\ref{wg44}) and (\ref{wg40}), we have 
the distributional inequality
$\alpha|\partial_t\chi|\le\liminf_{h\downarrow0}\frac{1}{\sqrt{h}}|\delta\chi|$,
provided the r.~h.~s.~is a finite measure.
Hence we obtain from (\ref{wg33}), (\ref{wg37}), and (\ref{wg55}) the distributional inequality
\begin{align}\label{wg58}
2\alpha\int_{0\le z\cdot\nu_0\le \alpha V_0}G_1dz|\partial_t\chi|
\le &\alpha^2\mu\notag\\
&+
2\int_{0\le z\cdot\nu_0\le \alpha V_0}G_1|z\cdot\nu|dz|\nabla\chi|dt,
\end{align}
which in particular shows that $\partial_t\chi$ is a measure.
Letting $V_0\uparrow\infty$ and appealing to (\ref{wg03}), 
this yields in particular $\alpha|\partial_t\chi|$ $\le\alpha^2\mu+4c_0|\nabla\chi|dt$
which we divide by $\alpha$:
\begin{align*}
|\partial_t\chi|\le\alpha\mu+\frac{4}{\alpha}c_0|\nabla\chi|dt.
\end{align*}
Letting $\alpha\downarrow 0$, we learn from the latter estimate
that null sets of $|\nabla\chi|dt$ are null sets of $\partial_t\chi$, so
that there exists $V\in L^1(|\nabla\chi|dt)$ such that (\ref{wg57}) holds.

\medskip

Since $|\partial_t\chi|=|V||\nabla\chi|dt$ is absolutely continuous
w.~r.~t.~to $|\nabla\chi|dt$, (\ref{wg58}) even holds with $\mu$ replaced
by its absolutely continuous part $\mu'$ w.~r.~t.~$|\nabla\chi|dt$.
Writing $\mu'=\theta|\nabla\chi|dt$ with $\theta\in L^1(|\nabla\chi|dt)$,
(\ref{wg58}) assumes the form
\begin{align*}
2\alpha\int_{0\le z\cdot\nu_0\le\alpha V_0}G_1dz|V|
\le\alpha^2\theta+
	2\int_{0\le z\cdot\nu_0\le\alpha V_0}G_1|z\cdot\nu|dz\quad|\nabla\chi|dt\mbox{-a.~e.}.
\end{align*}
As in Step \ref{Pr1S1}, a separability argument now allows to choose $\nu_0=\nu$,
so that by radial symmetry of $G_1$, the above assumes the form
\begin{align*}
2\alpha\int_{0\le z_1\le\alpha V_0}G_1dz|V|
\le\alpha^2\theta+
2\int_{0\le z_1\le\alpha V_0}G_1 z_1dz\quad|\nabla\chi|dt\mbox{-a.~e.}
\end{align*}
Dividing by $\alpha^2$ and momentarily writing $\alpha':=\alpha V_0$, this turns into
\begin{align*}
2\frac{V_0}{\alpha'}\int_{0\le z_1\le\alpha'}G_1dz|V|
\le\theta+2\frac{V_0^2}{{\alpha'}^2}\int_{0\le z_1\le\alpha'}G_1 z_1dz\quad|\nabla\chi|dt\mbox{-a.~e.}
\end{align*}
%
%
%
%
We now appeal to the limiting relations (which follow from factorizing $G_1$ into
the $(d-1)$-dimensional standard Gaussian and $G^{d=1}_1$)
\begin{align*}
	\lim_{\alpha'\downarrow 0}\frac{1}{\alpha'}\int_{0\le z_1\le\alpha'}G_1dz&=G_1^{d=1}(0)=c_0,\\
\lim_{\alpha'\downarrow 0}\frac{1}{{\alpha'}^2}\int_{0\le z_1\le\alpha'}G_1z_1dz
	&=\frac{1}{2}G_1^{d=1}(0)=\frac{c_0}{2},
\end{align*}
to see that the above turns into
\begin{align*}
	2c_0V_0|V|\le\theta+c_0V_0^2\quad|\nabla\chi|dt\mbox{-a.~e.}
\end{align*}
Again, by a separability argument for $V_0$, we may assume $V_0=|V|$
so that the above yields (\ref{wg53}) in form of $c_0V^2\le\theta$.

\bigskip

{\sc Proof of Proposition \ref{Pr3}}.

\newcounter{Pr3S} 
\refstepcounter{Pr3S} 

{\sc Step} \arabic{Pr3S}.\label{Pr3S1}\refstepcounter{Pr3S} Metric slope and functional derivative.
We claim the following relation between the metric slope $|\partial E(u)|$ of a functional
$E$ on ${\mathcal M}$, cf.~(\ref{wg92}), 
at a configuration $u$, and its first variation $\delta E(u).\xi$
in direction of a smooth vector field $\xi$:
\begin{align}\label{wg25}
\frac{1}{2}|\partial E(u)|^2\ge\delta E(u).\xi-\frac{1}{2}\big(\delta d(u,\cdot)(u).\xi\big)^2.
\end{align}
As usual, first variation $\delta$ is defined by considering the curve $s\mapsto u_s$ of
configurations characterized via the transport equation (to be interpreted distributionally
or solved explicitly with help of the flow map $\Phi_s$ via $u_s\circ\Phi_s=u$)
\begin{align}\label{wg27}
\frac{\partial u_s}{\partial s}+\xi\cdot\nabla u_s=0\quad\mbox{and}\quad u_{s=0}=u,
\end{align}
and setting
\begin{align}\label{wg28}
\delta E(u).\xi:=\frac{d}{ds}\Big|_{s=0}E(u_s)\quad\mbox{and}\quad
	\delta d(u,\cdot)(u).\xi:=\frac{d}{ds}\Big|_{s=0}d(u,u_s),
\end{align}
with the understanding that both derivatives exist (and define linear functionals in $\xi$
which is the case for $E=E_h$ and $d(u,\cdot)=d_h(u,\cdot)$, 
cf.~Steps \ref{Pr3S2} and \ref{Pr3S3}).
Inequality (\ref{wg25})
is then a direct consequence of the definition $|\partial E(u)|$, cf.~(\ref{wg26}),
which yields
\begin{align*}
	|\partial E(u)|&\ge\limsup_{s\downarrow 0}\frac{(E(u)-E(u_s))_+}{d(u,u_s)}\nonumber\\
	&\ge\frac{\lim_{s\downarrow 0}\frac{1}{s}(E(u_s)-E(u))}
	{\lim_{s\downarrow 0}\frac{1}{s}d(u,u_s)}
	=\frac{\delta E(u).\xi}{\delta d(u,\cdot)(u).\xi},
\end{align*}
and Young's inequality.

\medskip

{\sc Step} \arabic{Pr3S}.\label{Pr3S2}\refstepcounter{Pr3S}
Representation of $\delta E_h(u).\xi$; we claim:
%
\begin{align}\label{wg10}
\delta E_h(u).\xi=\frac{1}{\sqrt{h}}\int\Big(
	&\nabla\cdot\xi\big((1-u)G_h*u+uG_h*(1-u)\big)\nonumber\\
       +&u[\xi,\nabla G_h*](1-u)\Big),
\end{align}
where $[\xi,\nabla G_h*]$ $=\sum_{i=1}^d[\xi_i,\partial_iG_h*]$
denotes the commutator of multiplying with
$\xi$ and convolving with $\nabla G_h$. 
In checking this formula we may by approximation assume that $u$ is smooth;
by definition (\ref{wg27}) \& (\ref{wg28}) of $\delta$ we obtain from
the definition (\ref{wg06}) of $E_h$
\begin{align*}
\delta E_h(u).\xi=\frac{1}{\sqrt{h}}
\int\big((\xi\cdot\nabla u) G_h*u-(1-u)G_h*(\xi\cdot\nabla u)\big),
\end{align*}
which by the symmetry of $G_h*$ we rewrite as
\begin{align*}
\delta E_h(u).\xi=-\frac{1}{\sqrt{h}}
	\int\big(\xi\cdot\nabla(1-u)\,G_h*u+\xi\cdot\nabla u\,G_h*(1-u)\big).
\end{align*}
We write $\xi\cdot\nabla u=\nabla\cdot(u\xi)-u\nabla\cdot\xi$
and $\xi\cdot\nabla(1-u)=\nabla\cdot((1-u)\xi)-(1-u)\nabla\cdot\xi$,
so that (\ref{wg10}) reduces to the identity
\begin{align*}
-\int\big(\nabla\cdot((1-u)\xi)\,G_h*u&+\nabla\cdot(u\xi)\,G_h*(1-u)\big)\nonumber\\
&=\int u[\xi,\nabla G_h*](1-u),
\end{align*}
which follows from integration by parts and the anti-symmetry of $\nabla G_h*$.

\medskip

{\sc Step} \arabic{Pr3S}.\label{Pr3S3}\refstepcounter{Pr3S}
Representation of $\delta d_h(u,\cdot)(u).\xi$:
\begin{align}\label{wg19}
\lefteqn{\frac{1}{2}\big(\delta d_h(u,\cdot)(u).\xi\big)^2}\nonumber\\
&=\sqrt{h}\int\Big(u\xi\cdot\nabla^2G_h*((1-u)\xi)
	-u\xi\cdot\nabla G_h*((1-u)\nabla\cdot\xi)\nonumber\\
&+u\nabla\cdot\xi\,\nabla G_h*((1-u)\xi)-u\nabla\cdot\xi\,G_h*((1-u)\nabla\cdot\xi)\Big).
\end{align}
In the notation of Step \ref{Pr3S2},
$\frac{1}{2}\big(\delta d_h(u,\cdot)(u).\xi\big)^2 =\frac{1}{2}\big(\frac{d}{ds}\big|_{s=0}d_h(u,u_s)\big)^2$ 
$=\frac{1}{4}\frac{d^2}{ds^2}\big|_{s=0}d_h^2(u,u_s)$, so that by definition
(\ref{wg18}) of $d_h$ and by (\ref{wg27}), we have
$\frac{1}{2}\big(\delta d_h(u,\cdot)(u).\xi\big)^2$
$=\sqrt{h}\int\frac{\partial u_s}{\partial s}\big|_{s=0} \,G_h*\frac{\partial u_s}{\partial s}\big|_{s=0} $
$=\sqrt{h}\int\xi\cdot\nabla u\,G_h*(\xi\cdot\nabla u)$. Rewriting the second factor
$\xi\cdot\nabla u$
$=\nabla\cdot(u\xi)-u\,\nabla\cdot\xi$ and using the symmetry of $G_h*$,
an integration by parts, and $-\nabla u=\nabla(1-u)$, we obtain
\begin{align*}
\lefteqn{\frac{1}{2}\big(\delta d_h(u,\cdot)(u).\xi\big)^2}\nonumber\\
&=\sqrt{h}\int\big(u\xi\cdot\nabla G_h*(\xi\cdot\nabla(1-u))
	+u\nabla\cdot\xi\,G_h*(\xi\cdot\nabla(1-u))\big).
\end{align*}
We write $\xi\cdot\nabla(1-u)$ $=\nabla\cdot((1-u)\xi)-(1-u)\nabla\cdot\xi$
to obtain (\ref{wg19}).

\medskip

{\sc Step} \arabic{Pr3S}.\label{Pr3S3bis}\refstepcounter{Pr3S}
Passage to the limit in $\delta E_h$; we claim that
\begin{align}\label{wg29}
\lim_{h\downarrow 0}\int_0^T\delta E_h(u_h).\xi dt
=c_0\int(\nabla\cdot\xi-\nu\cdot\nabla\xi\nu)|\nabla\chi|dt.
\end{align}
According to (\ref{wg10}), we may split into two statements. The first statement is
\begin{align*}
\lim_{h\downarrow0}\frac{1}{\sqrt{h}}\int\nabla\cdot\xi\big((1-u_h)\,G_h*u_h
+u_h\,G_h*(1-u_h)\big)\nonumber\\
=2 c_0\int\nabla\cdot\xi|\nabla\chi|dt,
\end{align*}
which is an immediate consequence of testing (\ref{wg11}) \& (\ref{wg12})
with $\nabla\cdot\xi$, appealing to the scaling (\ref{wg05}) and
to the formula (\ref{wg03}). The second statement is
\begin{align}\label{wg13}
\lefteqn{\lim_{h\downarrow0}\frac{1}{\sqrt{h}}\int u_h[\xi,\nabla G_h*](1-u_h)}\nonumber\\
&=-c_0\int(\nu\cdot\nabla\xi\nu+\nabla\cdot\xi)|\nabla\chi|dt,
\end{align}
for which we now give the argument. Spelling out
\begin{align*}
\lefteqn{\big([\xi,\nabla G_h*](1-u_h)\big)(t,x)}\nonumber\\
&=\int(\xi(t,x)-\xi(t,x-z))\cdot\nabla G_h(z)(1-u_h)(t,x-z)dz,
\end{align*}
we see that
\begin{align}\label{wg15}
\big|&\big([\xi,\nabla G_h*](1-u_h)\big)(t,x)\nonumber\\
&-\int\nabla G_h(z)\cdot\nabla\xi(t,x) z(1-u_h)(t,x-z)dz\big|\nonumber\\
&\le\frac{1}{2}\sup|\nabla^2\xi| \int |z|^2|\nabla G_h(z)|(1-u_h)(t,x-z)dz.
\end{align}
Appealing to (\ref{wg05}) in form of $\nabla G_h(z)dz$ 
$=\frac{1}{\sqrt{h}}\nabla G_1(\frac{z}{\sqrt{h}})d\frac{z}{\sqrt{h}}$, we learn that the limit
of the contribution of the main term can be computed by
testing (\ref{wg11}) with $\frac{\nabla G_1(z)\cdot\nabla\xi(t,x) z}{G_1(z)}$
$=-z\cdot\nabla\xi(t,x)z$ (which is of polynomial growth in $z$); it assumes the value
\begin{align*}
\int\nabla G_1\cdot\nabla\xi z(\nu\cdot z)_+ |\nabla\chi|dtdz,
\end{align*}
which yields (\ref{wg13}) by formula (\ref{wg14}) below.
The contribution of the r.~h.~s.~error term in (\ref{wg15}) is vanishing of
$O(\sqrt{h})$, as follows from appealing to (\ref{wg11})
tested with $\frac{|\nabla G_1(z)||z|^2}{G_1(z)}$ $=|z|^3$.

\medskip

{\sc Step} \arabic{Pr3S}.\label{Pr3S5}\refstepcounter{Pr3S}
Passage to the limit in $\delta d_h$; we claim
\begin{align}\label{wg30}
\lim_{h\downarrow 0}\frac{1}{2}\int_0^T\big(\delta d_h(u_h,\cdot)(u_h).\xi\big)^2dt
	=c_0\int(\xi\cdot\nu)^2|\nabla\chi|dt.
\end{align}
According to (\ref{wg19}), this statement may be split into a leading-order statement
\begin{align}
\lim_{h\downarrow 0}\sqrt{h}\int u_h\xi\cdot\nabla^2G_h*((1-u_h)\xi) dxdt
&=c_0\int(\xi\cdot\nu)^2|\nabla\chi|dt,\label{wg20}
\end{align}
and the higher-order statements
\begin{align}\label{wg24}
\int u_h\xi\cdot\nabla G_h*((1-u_h)\nabla\cdot\xi)dxdt
	&=O(1),\\
\int u_h\nabla\cdot\xi\,\nabla G_h*((1-u_h)\xi) dxdt
	&=O(1),\nonumber\\
\frac{1}{\sqrt{h}}\int u_h\nabla\cdot\xi\,G_h*((1-u_h)\nabla\cdot\xi)dxdt
	&=O(1).\nonumber
\end{align}
The statement (\ref{wg20}) itself splits into the main part
\begin{align*}
\lim_{h\downarrow 0}\sqrt{h}\int \xi\cdot(u_h\nabla^2G_h*(1-u_h))\xi dxdt
=c_0\int(\xi\cdot\nu)^2|\nabla\chi|dt
\end{align*}
and the higher-order commutator
\begin{align}\label{wg23}
\int u_h\xi\cdot[\xi,\nabla^2G_h]*(1-u_h) dxdt=O(1).
\end{align}
The main part follows from appealing to (\ref{wg05})
in form of $\sqrt{h}\nabla^2G_h(z) dz$ 
$=\frac{1}{\sqrt{h}}\nabla^2 G_1(\frac{z}{\sqrt{h}}) d\frac{z}{\sqrt{h}}$,
testing (\ref{wg11}) with $\frac{\xi(t,x)\cdot\nabla^2 G_1(z)\xi(t,x)}{G_1(z)}$
$=(\xi(t,x)\cdot z)^2-|\xi(t,x)|^2$, 
and appealing to formula (\ref{wg21}).
The estimate of the error term (\ref{wg23}) follows from estimating the integrand by
$\sup|\xi|\sup|\nabla\xi|$ $\int |z||\nabla^2G_h(z)|\,u_h(t,x)(1-u_h)(t,x-z) dz$,
and then using the scaling (\ref{wg05}) 
further by 
\begin{align*}
\sup|\xi|\sup|\nabla\xi|\,\int |z||\nabla^2 G_1(z)|\,
\frac{1}{\sqrt{h}}u_h(t,x)(1-u_h)(t,x-\sqrt{h}z) dz,
\end{align*}
so that another application of (\ref{wg11}) yields (\ref{wg23}).
Statement (\ref{wg24}) and the other two higher-order estimates
follow along the same lines: For instance, the integrand in (\ref{wg24})
is $\le\sup|\xi|\sup|\nabla\cdot\xi|$ $\int|\nabla G_h(z)|$ $\,u_h(t,x)$ $(1-u_h)(t,x-z) dz$.
By rescaling, $\int|\nabla G_h(z)|$ $\,u_h(t,x)$ $(1-u_h)(t,x-z)$ $dz$
$=\int|z|G_1(z)\,\frac{1}{\sqrt{h}}u_h(t,x)(1-u_h)(t,x-\sqrt{h}z) dz$,
which is $O(1)$ by (\ref{wg11}).

\medskip

{\sc Step} \arabic{Pr3S}.\label{Pr3S6}\refstepcounter{Pr3S}
Conclusion. By Riesz' representation theorem in $L^2(|\nabla\chi|dt)$ and an approximation argument 
in the (arbitrary) smooth vector field $\xi$, the statement 
of Proposition \ref{Pr3} is a consequence of
\begin{align*}
\liminf_{h\downarrow0}\int_0^T\frac{1}{2}|\partial E_h(u_h)|^2dt\ge
c_0\int\big(\nabla\cdot\xi-\nu\cdot\nabla\xi\nu-(\xi\cdot\nu)^2\big)|\nabla\chi|dt.
\end{align*}
The latter follows starting from the inequality (\ref{wg25}) for $(E,d,u)=(E_h,d_h,u_h(t))$,
integrating in $t\in(0,T)$, and appealing to (\ref{wg29}) and (\ref{wg30}) to
pass to the limit on the r.~h.~s.~.

\medskip

{\sc Step} \arabic{Pr3S}.\label{Pr3S8}\refstepcounter{Pr3S}
Two formulas: For any unit vector $\nu$, any matrix $A$, and any vector $\xi$, we have
\begin{align}
-\int\nabla G_1(z)\cdot A z\,(\nu\cdot z)_+dz&=c_0(\nu\cdot A\nu+{\rm tr A}),\label{wg14}\\
\int\xi\cdot\nabla^2G_1(z)\xi\,(\nu\cdot z)_+dz&=c_0(\xi\cdot\nu)^2\label{wg21}.
\end{align}
Since $G_1(z)=-zG_1(z)$, for the first formula
we may assume that $A$ is symmetric; by linearity
we may assume that $A=e\otimes e$ for some unit vector $e$; by radial symmetry of $G_1$,
it thus remains to show
\begin{align*}
-\int\partial_1 G_1 z_1(\nu\cdot z)_+dz=c_0(\nu_1^2+1),
\end{align*}
which by one integration by parts, taking into account (\ref{wg03}), reduces to
\begin{align*}
\int_{\nu\cdot z>0} G_1 z_1dz=c_0\nu_1,
\end{align*}
which in view of $G_1 z_1=-\partial_1G_1=-\nabla\cdot(G_1 e_1)$ by the divergence theorem 
reduces to 
\begin{align}\label{wg22}
\int_{\nu\cdot z=0} G_1=c_0,
\end{align}
which follows since $c_0=G_1^{d=1}(0)$.
We now turn to (\ref{wg21}). By radial symmetry of $G_1$ and homogeneity, it suffices to
show
\begin{align}
\int\partial_1^2G_1(\nu\cdot z)_+dz=c_0\nu_1^2,
\end{align}
which by two integration by parts reduces to (\ref{wg22}).

\bibliographystyle{plain}

\bibliography{lit}

\begin{thebibliography}{10}

\bibitem{AlmgrenTaylorWang}
Fred Almgren, Jean~E. Taylor, and Lihe Wang.
\newblock Curvature-driven flows: a variational approach.
\newblock {\em SIAM Journal on Control and Optimization}, 31(2):387--438, 1993.

\bibitem{AmbrosioGigliSavare}
Luigi Ambrosio, Nicola Gigli, and Giuseppe Savar{\'e}.
\newblock {\em Gradient flows in metric spaces and in the space of probability
  measures}.
\newblock Birkh\"auser, 2008.

\bibitem{EsedogluOtto}
Selim Esedo\u{g}lu and Felix Otto.
\newblock Threshold dynamics for networks with arbitrary surface tensions.
\newblock {\em Communications on Pure and Applied Mathematics}, 68(5):808--864,
  2015.

\bibitem{Evans}
Lawrence~C. Evans.
\newblock Convergence of an algorithm for mean curvature motion.
\newblock {\em Indiana University Mathematics Journal}, 42(2):533--557, 1993.

\bibitem{LauxOtto1}
Tim Laux and Felix Otto.
\newblock Convergence of the thresholding scheme for multi-phase
  mean-cur\-va\-ture flow.
\newblock {\em Calculus of Variations and Partial Differential Equations},
  55(5):1--74, 2016.

\bibitem{LauxOtto2}
Tim Laux and Felix Otto.
\newblock Brakke's inequality for the thresholding scheme.
\newblock {\em arXiv:1708.03071 [math.AP]}, 2017.

\bibitem{LuckhausSturzenhecker}
Stephan Luckhaus and Thomas Sturzenhecker.
\newblock Implicit time discretization for the mean curvature flow equation.
\newblock {\em Calculus of Variations and Partial Differential Equations},
  3(2):253--271, 1995.

\bibitem{BenceMerrimanOsher}
Barry Merriman, James~K. Bence, and Stanley~J. Osher.
\newblock Diffusion generated motion by mean curvature.
\newblock CAM Report 92-18, 1992.
\newblock Department of Mathematics, University of California, Los Angeles.

\bibitem{MichorMumford}
Peter~W. Michor and David Mumford.
\newblock Riemannian geometries on spaces of plane curves.
\newblock {\em Journal of the European Mathematical Society}, 8(1):1--48, 2006.

\bibitem{EsedogluSalvador}
Tiago Salvador and Selim Esedo\u{g}lu.
\newblock A simplified threshold dynamics algorithm for isotropic surface
  energies.
\newblock {\em Journal of Scientific Computing}, pages 1--22, 2018.

\end{thebibliography}

\end{document}